\documentclass[12pt,a4paper]{amsart}
\usepackage{amsmath, amssymb, bm}
\usepackage[all]{xy}

\def \LL {\skp \setlength{\leftskip}{8pt}}  
\def \LB {\setlength{\leftskip}{0pt}}   
\let \lan \langle
\let \ran \rangle

\def \c {\colon}
\def \cc {\, \colon \!}

\def \q {\qquad}
\def \qq {\qquad \qquad}
\def \qqq {\qquad \qquad \qquad \qquad}

\def \sep {\,  |  \,}     
\def \bu {{\scriptscriptstyle\bullet}}
\def \skp {\medskip}
\def \ndt {\noindent}
\def \Ndt {\medskip  \noindent}
\def\shp {^{\sharp}}
\def \tilde {\raise.17ex\hbox{$\scriptstyle\mathtt{\sim}$}}   

\def \rlw {\;\; \raisebox{0.4ex}{$\longrightarrow$} \hspace{-3.7ex}   
 \raisebox{-0.4ex}{$\longleftarrow$} \;\;}

\def \ard {\ar@{-->}}
\def \arp {\ar@{.>}}
\def \arv {\ar@{}}   
\def \arl {\ar@{-}}    
\def \arld {\ar@{--}}    
\def \arlp {\ar@{..}}    

\def \sup {{\rm sup}}
\def \inf {{\rm inf}}
\def \ti {\! \times \!}
\def \te {\otimes}
\def \jo {{\, {\scriptstyle{\vee}} \,}}
\def \Jo {{\, {\vee} \,}}
\def \me {{\, {\scriptstyle{\wedge}}\, }}
\def \Me {{\, {\wedge}\, }}
\def \ci {{\raise.3ex\hbox{$  \scriptscriptstyle\circ $}}}
\def \Pro {\raisebox{0.45ex}{${\mbox{\fontsize{10}{10}\selectfont\ensuremath{\prod}}}$}}
\def \Sum {\raisebox{0.45ex}{${\mbox{\fontsize{10}{10}\selectfont\ensuremath{\sum}}}$}}

\def \Ten {\raisebox{0.42ex}{${\mbox{\fontsize{10}{10}\selectfont\ensuremath{\bigotimes}}}$}}
\def \setm {{\raise.4ex\hbox{$ \, \scriptscriptstyle{\setminus} \; $}}}
\def \sub {\subset}
\def\le{\leqslant}
\def\ge{\geqslant}
\def \iso {\: \cong \:}

\def \and {\mbox{ and }}
\def \for {\mbox{for }}

\def \inn {\mbox{ in }}
\def \IF {\; \mbox{ if } \:}

\def \ga {\gamma}
\def \de {\delta}
\def \ep {\varepsilon}

\def \la {\lambda}
\def \si {\sigma}
\def \ph {\varphi}
\def \De {\Delta}

\def \id {{\rm id\,}}
\def \op {^{{\rm op}}}
\def \tp {{\rm tp}}
\def \mt {{\rm mt}}
\def \pr {{\rm pr}}
\def \Up {{\underline{p}}} 
\def \Ut {{\underline{t}}}
\def \Ux {{\underline{x}}}
\def \Uy {{\underline{y}}}
\def \Uz {{\underline{z}}}
\def \Oinfty {{\overline{\infty}}}
\def \Umu {{\underline{\mu}}}

\def \Set {\mathsf{Set}}
\def \Top {\mathsf{Top}}
\def \Mtr {\mathsf{Mtr}}

\def \pSet {{\rm p}\mathsf{Set}}
\def \wPi {{\rm w} \Pi}
\def \bbX {\mathbb{X}}

\def \es {\varnothing}
\def \sing  {\{*\}}

\def \bbZ {\mathbb{Z}}
\def \bbI {\mathbb{I}}
\def \bbR {\mathbb{R}}
\def \bbS {\mathbb{S}}
\def \OR {\overline{\mathbb{R}}}
\def \Bw {{\bf w}}
\def \vR {{\rm v} \mathbb{R}}
\def \vS {{\rm v} \mathbb{S}}

\begin{document}

\title[Weighted algebraic topology, II] {Weighted algebraic topology, II \\ 
(Real valued metrics)}

\author[M. Grandis]{Marco Grandis}

\address{Marco Grandis, Dipartimento di Matematica, Universit\`a di Genova, 16146-Genova, Italy.}
\email{grandismrc@gmail.com}



\begin{abstract}

Extending the `metric spaces' of Lawvere, we study `real metrics', with values in the extended 
real line $[- \infty, \infty]$.

	Formally, this ordered set is a symmetric monoidal closed category, and 
our structures are enriched categories on the latter. Concretely, the present goal is measuring 
`profits' and `losses' of a process, in any sense -- possibly related to energy, or a variable in 
any science. In particular, {\em linear} real metrics derive from a potential function. 

	This article is Part II in a series devoted to `weighted algebraic topology' -- an enriched 
version of directed algebraic topology, where paths are measured. Part III will introduce a finer
framework, more adequate to `quotient spaces' (as the spheres) and better related to topology.

\end{abstract}

 \maketitle

\section*{Introduction}\label{Intro}

\subsection{Metrics with real values}\label{0.1}
	This article is Part II of a series which began with \cite{G3}, in 2007, although at the time 
a sequel was not foreseen. The matter of the latter, reorganised in the general context of 
Directed Algebraic Topology, can also be found in Chapter 6 of the book \cite{G4}. Both texts 
are freely available.

	A prime structure analysed in \cite{G3} consists of Lawvere's generalised metrics \cite{L1}, 
with values in the extended positive line $\OR^+ = [0, \infty]$. Here we further extend this 
notion, studying `real metrics', with values in the extended real line $ \OR = [- \infty, \infty]$, 
to express `profits' and `losses', in any sense -- possibly related to a variable of any science,
from Physics to Economy.

	Two unpublished preprints of 1984 and 2015 also use metrics with values in $\OR$, 
in different perspectives. Dealing with entropy, Lawvere introduced a real valued 
antimetric $M(x, y)$ on the class of objects of a category $\bbX$ endowed with an 
`entropy supply' $s\c \bbX \to \bbR$ (\cite{L2}, p. 4). Dealing with the Legendre-Fenchel 
transform for real vector spaces, Willerton used real valued metrics, called 
$\OR$-{\em metrics} \cite{Wi}.

\subsection{An elementary example}\label{0.2}
	The present goals should be evident from the following elementary example, which does 
not require any prerequisite -- provided we are willing to follow its concrete interpretation; 
definitions and results on this matter will be made precise in Section 5, after studying the 
general framework.

\Ndt (a) We consider the vertical profile $a$ of a track, or road (in a limited area of the Earth),
represented in the $xy$-plane, with points $p' = (x', y')$, $p'' = (x'', y'')$, $p_j = (x_j, y_j)$
%
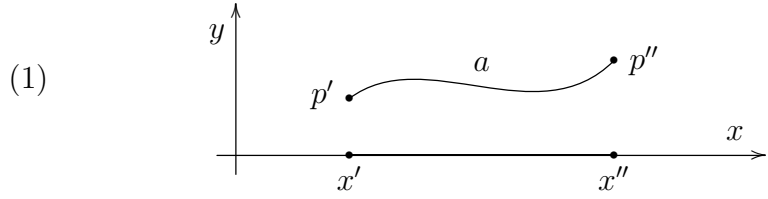
    \begin{equation} 
\xy <.5mm, 0mm>:
(3,-4) *{p'}; (10,-26) *{x'}; (88,5) *{p''}; (80,-26) *{x''}; (45,4) *{a}; (112,-14) *{x}; (-25,12) *{y}; 
(0,28) *{}; (0,-32) *{}; 
@i @={(10, -5), (10, -20), (80, 5), (80, -20)} @@{*{\bu}};
(10, -5); (80, 5) **\crv{(30,10)&(60,-15)},
\POS(-25,-20) \ar+(145,0),  \POS(-20,-25) \ar+(0,45), \POS(10,-19.8) \arl+(70,0),  
\endxy
    \label{0.2.1} \end{equation}

	As usual in such maps (for mountain paths, etc.), the horizontal projection on the 
$x$-axis represents the `rectified horizontal path', that is, the rectification of the projection 
of the track on a horizontal plane, where the distance is measured by 
$ \de_1(p_1, p_2) = |x_2 - x_1|$. The vertical axis measures the altitude, with distance 
$\, \de_2(p, q)$ = $|y_2 - y_1|$.

	The distance on the path itself, and its length, are measured by the euclidean metric of 
the plane $xy$
    \begin{equation} \begin{array}{c}
d(p, q)  =  (\de_1(p, q)^2 + \de_2(p, q)^2)^{1/2}
\\[6pt]
L_d(a)  =  \sup_\Up \, (\Sum_j \; d(p_{j-1}, p_j)),
    \label{0.2.2} \end{array} \end{equation}
by a least upper bound on all finite sequences $\Up = (p_0, ..., p_n)$ of points of the path 
{\em with increasing first coordinate}.

\Ndt (b) We are also interested in distinguishing ascent and descent -- supposing we are going 
from $ p' $ to $ p''$. This is not measured by the pseudometric $\de_2$, but can be measured 
by a a real valued metric $\rho$, and its positive part $\rho^+ = \rho \jo 0$
    \begin{equation}
\rho(p_1, p_2)  =  y_2 - y_1,   \q\;\;   \rho^+(p_1, p_2)  =  \max(y_2 - y_1, 0)  \ge  0.
    \label{0.2.3} \end{equation}

	Let us note that $ \rho $ is `linear' (it is finite and `additive', see \ref{2.2}(c)), 
and expresses the increment of the gravitational potential per mass-unit, while $ \rho^+ $ is a 
quasi-pseudo metric \cite{Ke} and a generalised metric in Lawvere's sense \cite{L1} 
(see \ref{0.3}). Both $ \rho $ and $ \rho^+ $ are invariant up to translation and not symmetric.

\Ndt (c) The real metric $ \rho$ gives a real valued {\em valuation} $ v_\rho(a) $ of the path, 
which, by linearity, is simply the increment in altitude from $ p' $ to $ p''$
    \begin{equation}
v_\rho(a)  =  \sup_\Up \, (\Sum_j \; \rho(p_{j-1}, p_j))  =  y'' - y'  \in  \bbR,
    \label{0.2.4} \end{equation}
with sequences $\Up$ varying as above.

	The {\em total ascent} $v^+_\rho(a)$ of the path is defined as its $\rho^+$-{\em valuation} 
    \begin{equation} \begin{array}{c}
v^+_\rho(a) =  v_{\rho^+}(a)  =  \sup_p \, (\Sum_j \; \rho^+(p_{j-1}, p_j)),
\\[5pt]
v^+_\rho(a) \, \in \, \bbR^+ = [0, \infty[,   \q   v^+_\rho(a)  \ge  v_\rho(a).
    \label{0.2.5} \end{array} \end{equation}

	We assume that the track is the graph of a piecewise monotone continuous function 
$a\c [x', x''] \to \bbR$. Thus $ v^+_\rho(a) $ is finite: the sum of increments in altitude on the 
maximal intervals on which the function $ a $ is increasing.

	The {\em total descent} $ v^-_\rho(a)$, defined as
    \begin{equation}
v^-_\rho(a)  =  v^+_\rho(a) - v_\rho(a)  \ge  0   \q\;\;   (v_\rho(a)  =  v^+_\rho(a) - v^-_\rho(a)),
    \label{0.2.6} \end{equation}
is also a real number. Owing again to the linearity of $ \rho$, $v^-_\rho(a) $ can also be obtained 
as the total ascent of the reversed path $a\shp$, going backwards (as proved in Lemma 
\ref{5.3})
    \begin{equation}
v^+_\rho(a\shp)  =  \sup_p \, (\Sum_j \; \rho^+(p_j, p_{j-1})).
    \label{0.2.7} \end{equation}

	The absolute value $ |\rho|(p_1, p_2) = |y_2 - x_2| $ is here a pseudometric (another 
consequence of the linearity of $ \rho$). The $|\rho|$-valuation of the path expresses the 
{\em total variation} of altitude
    \begin{equation}
v_{|\rho|}(a)  =  \sup_p(\Sum_j |\rho|(p_{j-1}, p_j)),   \q   v_{|\rho|}(a)  =  v^+_\rho(a) + v^-_\rho(a).
    \label{0.2.8} \end{equation}

\Ndt (d) Figure \eqref{0.2.1} might also represent the value $ y $ of any variable (in Physics, 
Chemistry, Natural Sciences, Economy, Sociology, etc.) as a function of the variable $ x $ 
(position on a line, or time, etc.). The interpretation of $ v_\rho(a)$, $ v^+_\rho(a)$, 
$ v^-_\rho(a)$, $ v_{|\rho|}(a) $ is similar.

\skp	More complex examples, related to the gravitational or electric potential in a cartesian 
space, will be examined in Section 6.

\subsection{Positive metrics and weighted algebraic topology}\label{0.3}
Lawvere's article \cite{L1} presented, in 1974, a (generalised) metric space $ X $ as a small 
category enriched on the positive extended real line $\Bw^+ = \OR^+ = [0, \infty]$; the latter is 
viewed as a category (with morphisms $\la \ge \rho$) endowed with a strict symmetric monoidal 
closed structure given by the (extended) sum $\la + \rho$.

	Concretely, $ X $ is a set equipped with a (generalised) metric $\de\c X \ti X \to \Bw^+$ 
(the enriched hom) satisfying the triangular inequality (enriched composition) and a second 
condition (enriched units)
    \begin{equation}
\de(x, y) + \de(y, z)  \ge  \de(x, z),   \q\;\;   0  \ge  \de(x, x)   \q   (x, y, z \in X),
    \label{0.3.1} \end{equation}
that amounts to $ \de(x, x) = 0$, for all $ x \in X$. Symmetry and `separation' are not required.

	Similar extensions have been studied in the theory of metric spaces, like the
{\em quasi-pseudo metrics} of \cite{Ke}, which take finite values, or the {\em \'ecarts} of 
Bourbaki \cite{Bk}, which are symmetric.

	Lawvere's extension combines a formal elegance within category theory (one might say 
`necessity') with concrete advantages, e.g.\ with respect to the existence of limits and colimits 
-- which obviously fails for `finite metrics'. (An influential part of the article \cite{L1}, developed 
by several authors, extends part of the theory of metric spaces, like Cauchy completeness, to 
enriched categories.)

	Breaking the symmetry property, such a metric can also express privileged directions. It 
was used in this sense by the present author, under the name of $\de$-{\em metric}, as a base 
for an enriched version of Directed Algebraic Topology \cite{G3, G4}. Lipschitz paths and 
Lipschitz homotopies are based on the {\em standard} $\de$-{\em interval} $\de\bbI = \de[0, 1]$, 
with a $\de$-metric meant to hinder any backward move
    \begin{equation}
\de(t, t')  =  t' - t  \IF  t \le t',   \qq   \de(t, t')  =  \infty  \; \text{ otherwise.}
    \label{0.3.2} \end{equation}

	The length of a Lipschitz path $ a\c \de\bbI \to X $ is given by the following least upper bound 
in $\OR^+$:
    \begin{equation} \begin{array}{c}
L(a)  =  \sup_\Ut \, L_\Ut (a),
\\[5pt]
L_\Ut(a)  =  \Sum_j \; \de(a(t_{j-1}), a(t_j))   \q\;\;   (0 = t_0 \le t_1 \le ... \le t_n =1),
    \label{0.3.3} \end{array} \end{equation}
where $\Ut = (t_j)$ stands for any finite increasing sequence in $[0, 1]$.

	A $\de$-metric space $X$, or $\de$-space, has a {\em fundamental weighted category} 
$\wPi_1(X) $ \cite{G3, G4}, an enrichment of the fundamental category of a directed spaces 
\cite{G1, G4}: it 
has for objects the points of $ X$, and for arrows $\xi\c x \to x' $ the classes of Lipschitz paths 
$ a\c \de\bbI \to X $ from $ x $ to $ x'$, up to the equivalence relation generated by Lipschitz 
homotopies with fixed endpoints. The {\em weight} of such an arrow $\xi\c x \to x' $ is defined 
as usual in a quotient, as the greatest lower bound of the length of the paths belonging to this 
equivalence class
    \begin{equation}
w(\xi)  =  \inf_{a\in \xi} \, L(a).
    \label{0.3.4} \end{equation}

	This gives $\wPi_1(X) $ a {\em weight}, or {\em generalised norm}, in $ \Bw^+ $ 
\cite{L1, G2, G3, G4}. All this is extended to the more general framework of `spaces with 
weighted paths', whose finer quotients can interpret structures of Noncommutative Geometry: 
see \cite{G3}, Section 7, and \cite{G4}, Section 6.7.

\subsection{Real valued metrics}\label{0.4}
We now replace $ \Bw^+ $ by a larger structure $ \Bw $ based on the extended real line 
$ \OR = [- \infty, \infty]$. This complete lattice is again viewed as a category, with morphisms 
$\la \ge \rho$, and endowed with a strict symmetric monoidal closed structure given by the 
extended sum $\la + \rho$; here we let $ - \infty + \infty = \infty $ (adding any element must 
preserve the {\em initial} object $ \infty$).

	A {\em $\rho$-metric space}, or {\em $\rho$-space}, is thus a small category $ X $ 
enriched on $ \Bw$. This means a set $ X $ with a mapping $ \rho\c X \ti X \to \OR$, called a 
{\em real metric}, or {\em $\rho$-metric}, that satisfies the same conditions as $\de$ in 
\eqref{0.3.1}
    \begin{equation}
\rho(x, y) + \rho(y, z)  \ge  \rho(x, z),   \q\;\;   0  \ge  \rho(x, x)   \q   (x, y, z \in X),
    \label{0.4.1} \end{equation}
but the second condition amounts here to $ \rho(x, x) = 0 $ {\em or} $ - \infty$, for all $ x \in X $ 
(see \ref{2.1}).

	The extended real number $ \rho(x, y) $ expresses a measure of the transition from 
$x$ to $y$. The `valuation', or `real-valued length', of a path is also an extended real number 
(see \ref{5.2}).

	The basic example is now the {\em standard} $\rho$-{\em line} $ \rho\bbR$, with 
real metric
    \begin{equation}
\rho(x, y)  =  y - x.
    \label{0.4.2} \end{equation}

	More generally, a {\em linear} $\rho$-metric on a set $ X $ (see \ref{2.2}(c)) can be 
expressed as the variation $ \rho(x, y) = \Phi(y) - \Phi(x) $ of an 
arbitrary function $ \Phi\c X \to \bbR$, called the {\em potential} of $ \rho $ (see \ref{3.2}).

	Linear real metrics are important here (and trivial in the positive case); they could also 
be called `conservative metrics', being related to a potential. These two features will be 
distinguished in Part III: the more general structure of a `space with valued paths' presents a 
{\em linear} case (possibly determined by a closed 1-form on a smooth manifold), which 
only corresponds to a linear $\rho$-metric in the {\em conservative} case (for an exact 1-form, 
or a conservative vector field on a Riemannian manifold).

\subsection{Outline}\label{0.5}
The enrichments we are interested in, over $ \Bw^+ $ and $ \Bw$, are made clear considering 
the structure of a {\em commutative quantale}, in Section 1. This can be equivalently defined as 
a complete lattice suitably enriched, or a complete, symmetric monoidal closed category of 
order type.

	Section 2 introduces $\rho$-metric spaces, developing their links to the positive 
case of $\de$-metric spaces previously studied in \cite{L1, G3, G4} and other works. Limits and 
colimits are dealt with in \ref{2.3}. The main properties, like linearity and invariance up to group 
actions, are outlined in \ref{2.2} and studied in Section 3. As a critical problem, quotients of 
$\rho$-spaces tend to be trivial: the line $ \rho\bbR $ induces a chaotic $\rho$-metric on the 
circle (with constant value $ - \infty$), see \ref{3.1}(c).

	Section 4 deals with the reflective and coreflective symmetrisations of a $\rho$-space, 
the associated topologies, and the associated reflective and coreflective preorders. In different 
situations it is convenient to investigate spaces and maps by `reflective' or `coreflective tools' 
-- another critical aspect discussed below. (In Part I only the reflective procedures were needed, 
for $\de$-spaces.)

	Paths in a $\rho$-space and their {\em valuation} in $\OR$ are analysed in Section 5, 
extending the {\em length} of paths in a $\de$-space. The tensor product of $\rho$-spaces is 
studied in Section 6, together with various examples, including real metrics associated to a 
conservative vector field.

\subsection{Comments and sequels}\label{0.6}
(a) Directed Algebraic Topology \cite{G4} studies `directed spaces' with `directed algebraic 
structures' produced by homotopy or homology functors: on the one hand the fundamental 
{\em category} (and its higher dimensional versions), on the other hand {\em preordered} 
homology groups. Its general aim is {\em modelling non-reversible phenomena}. It is currently 
used in the theory of concurrent systems, cf.\ \cite{FGHMR}.

	Weighted Algebraic Topology, as proposed in \cite{G3, G4}, is an enrichment of this topic: 
the truth-valued approach of directed algebraic topology (where a path is either allowed or not) 
is replaced with a measure of costs, possibly infinite. The general aim is, now, {\em measuring 
the cost of (possibly non-reversible) phenomena}. 

	Here we further extend this approach, so that a cost can also be negative, expressing a 
gain. Various aspects of the extended theory are drastically different from the positive case, and 
present critical problems already mentioned in the outline, with respect to quotients and topology.

\Ndt (b) The second point is the more intriguing: real metrics do have a topological content, 
studied throughout this paper, {\em but a complex relation to topology}. We should not expect a 
real metric space to have a {\em definite} topology: there is rather a range of them, from the 
reflective to the coreflective topology (see Section 4), and even this range can be insufficient, 
as in \ref{4.8}.  In our perspective, the main interest of real metrics is evaluating paths 
(see Section 5), but different goals have already appeared in \cite{L2, Wi}.

	A convincing evidence in this sense will be explored in a sequel. The Minkowski 
spacetime of special relativity can be given a real valued metric, derived from the metric tensor, 
much in the same way as a Riemannian manifold gets an ordinary metric. The path of a 
material or lightlike particle is continuous for the natural (euclidean!) topology, which has little 
to do with the metric tensor, and the associated real metric. On the other hand, these structures 
do evaluate the proper duration of a causal path in the spacetime.

\Ndt (c) Part III will overcome these problems, and others, introducing the more flexible setting 
of {\em valued spaces}, or {\em spaces with valued paths} -- a real valued version of the positive 
case already studied in \cite{G3, G4}. In this framework an object already bears a given 
topology, and quotients work as one can expect: for instance, the valued line $\vR$ (based 
on the present $\rho$-line $\rho\bbR$) has an appropriate quotient up to integral translations: 
the (linearly) valued circle $\vS^1$, where a path is valued by any lifting to the line $\rho\bbR$, 
and a loop is valued by its (integral) winding number; this also shows that, for valued spaces, 
the linear case need not be conservative, as already mentioned in \ref{0.4}.

	Homotopy theory for $\rho$-spaces and valued spaces will also be studied.

\subsection{Notation and terminology}\label{0.7}
Open and semiopen intervals of the real line are always denoted by square brackets, like 
$ ]0, 1[$, $[0, 1[$, etc. The symbol $\sub $ denotes weak inclusion.

	In a category, the term `map' is used as a synonym of `morphism'; in particular, between 
topological spaces, it is always used in the sense of {\em continuous mapping}. In an ordered set, 
joins and meets are written as $ x \jo y $ and $ x \me y $ (also for the minimum and maximum of a
pair, in a total ordering). Increasing and decreasing functions are always meant in the weak 
sense.

	The `extended numbers' $ \pm \infty $ are written $ \infty $ and $ \Oinfty$, to avoid ambiguous 
expressions (see \ref{1.1}(c)).

	Marginal remarks and elementary proofs are written in small characters. For the sake of 
brevity, some technical points are referred to the preprint version of this article \cite{G5}.

\subsection{Acknowledgements}\label{0.8}
I am indebted to my colleagues Ettore Carletti and Fernando Gliozzi, for long discussions on the present subject.

\section{Values and weights}\label{s1}
	We begin by examining the extended real line $ \OR = [\Oinfty, \infty] $ and the extended 
half-line $ \OR^+ = [0, \infty]$, containing the values of the real or positive metrics dealt with 
in the next section. The terms `real' and `positive' will generally refer to these ordered sets.

	As in Lawvere's article \cite{L1}, based on the extended half-line, lattice terminology will 
refer to the natural order of $ \OR $ and $ \OR^+$, while categorical terms will refer to the direction 
$\la \ge \mu$ of morphisms. This (unfortunate!) reversion is required by using numbers as values 
of metrics; for instance, the inequality $ d \ge d'$ between two metrics on the same set 
$X$ gives a `comparison map' $ (X, d) \to (X, d') $ that lifts the identity. (Dealing with antimetrics 
would overcome this issue, but complicate the relationship with classical metrics.)

\subsection{Values, or real weights}\label{1.1}
(a) The real metrics we shall introduce take values in the extended line $\OR = [\Oinfty, \infty]$, 
or -- more precisely -- in a structure $\Bw$ which can be analysed in two ways, based either on 
lattices and quantales \cite{Mu, Ro, Ye} or on monoidal closed categories \cite{EK, Ma, Ky}. A 
reader can focus on one of these aspects, but both are important: many points are more concrete 
in the terminology of lattices, but better motivated by category theory. (Some verifications are 
deferred to point (b).)

\LL

\ndt (i) As an ordered set, $\Bw$ is the extended line $\OR = [\Oinfty, \infty]$ with the natural order 
$ \la \le \mu $ and the {\em extended sum} $ \la + \mu$, where $ \Oinfty + \infty = \infty$. It is a 
commutative unital coquantale, that is, a complete lattice with an operation $ \la + \mu $ of 
abelian monoid which preserves arbitrary meets in each variable (and, in particular, the 
maximum $\infty$).

\Ndt (ii) As a category, $ \Bw $ has objects $ \la \in \OR$, a unique morphism $ \la \to \mu $ when 
$ \la \ge \mu $ and a monoidal structure given by the extended sum $ \la + \mu$. It is a small 
complete category of order type, with a symmetric monoidal closed structure.

 \LB

\begin{small}

\Ndt (b) To show that $ \Bw $ is indeed an ordered abelian monoid (and a symmetric monoidal 
category) it is sufficient to note that the additive group $ \bbR $ (or any ordered abelian monoid) 
can be augmented adding an absorbent minimum and {\em then} an absorbent maximum.

	Moreover, the fact that each mapping $ \la + (\; )\c \Bw \to \Bw $ preserves all meets 
(or colimits) is equivalent to the existence of the internal hom (computed below), by the 
Adjoint Functor Theorem.

\end{small}

\Ndt (c) Joins and meets in $\Bw$ give `all' limits and colimits in our category, which amount to 
products and sums
    \begin{equation} \begin{array}{llll}
\text{product:}   &   \Jo \la_i,   \q   &   \text{terminal object:}   &  \Oinfty,
\\[5pt]
\text{sum:}   &  \Me \la_i,   &   \text{initial object:}   &   \infty.
    \label{1.1.1} \end{array} \end{equation}

	The internal hom-object $\, \hom(\mu, \nu) \, $ is given by an {\em extended difference}, 
also written as $\, \nu - \mu$, and characterised as follows
    \begin{equation} \begin{array}{c}
\la + \mu \ge \nu   \;  \iff \;    \la \ge  \hom(\mu, \nu),
\\[5pt]
\hom(\mu, \nu)  =  \nu - \mu  =  \min\{\la \in [\Oinfty, \infty] \, \sep \la + \mu \ge \nu\}.
    \label{1.1.2} \end{array} \end{equation}

	The following `indeterminate forms' are thus determined:
    \begin{equation} \begin{array}{c}
\infty - \infty  =  \min\{\la \in [\Oinfty, \infty] \, \sep \la + \infty \ge \infty\}  =  \Oinfty,
\\[5pt]
\Oinfty - \Oinfty  =  \min\{\la \in [\Oinfty, \infty] \, \sep \la + \Oinfty \ge \Oinfty\}  =  \Oinfty.
    \label{1.1.3} \end{array} \end{equation}

	The notation $ \Oinfty $ allows us to distinguish the extended sum $\infty + \Oinfty = \infty$ 
and the extended difference $ \infty - \infty = \Oinfty$. Note also that $\infty$ has no additive 
inverse in $\Bw$.

\begin{small}

\Ndt (d) The category $ \Bw $ is *-autonomous \cite{Ba}, with dualising object the additive 
identity 0 and involution $ \la^* = \hom(\la, 0) = - \la$.

\end{small}

\subsection{Positive values, or weights}\label{1.2}
(a) The complete lattice $ \Bw^+ = \OR^+ = [0, \infty] $ of extended positive numbers, also called 
{\em positive values}, or {\em weights}, inherits a monoidal structure $ \la + \mu $ which is also 
closed (in its own right). It was introduced by Lawvere [L1], in 1974 (the original notation is 
${\bf R}$). Weight 0 should be viewed as `free', and weight $ \infty $ as `unreachable, 
unaffordable, or forbidden'. Metrics with values in $ \Bw^+ $ will be reviewed in \ref{2.4}.

	The internal hom, written here as $\, \hom^+(\mu, \nu)$, is given by the {\em truncated 
difference} $ (\nu - \mu) \jo 0$
    \begin{equation} \begin{array}{c}
\la + \mu \ge \nu   \; \iff \;   \la \ge  \hom^+(\mu, \nu),
\\[5pt]
\hom^+(\mu, \nu)  =  \min\{\la \in [0, \infty] \sep \la + \mu \ge \nu\}  =  (\nu - \mu) \jo 0,
\\[5pt]
\hom^+(\infty, \infty)  =  \min\{\la \in [0, \infty] \, \sep \la + \infty \ge \infty\}  =  0.
    \label{1.2.1} \end{array} \end{equation}

\begin{small}

	Truncated difference is written $ \nu - \mu $ in \cite{L1}, a notation introduced above 
in a different sense.

\end{small}

\Ndt (b) The embedding $ \Bw^+ \sub \Bw $ is a strictly monoidal functor: it preserves the
tensor product, meets (i.e.\ colimits) and non-empty joins, but not the minimum (the terminal 
object). It has a coreflector (right adjoint), called {\em positive part} (of a value)
    \begin{equation}
( \;\; )^+\c \Bw \to \Bw^+,   \q\;\;   \la^+  =  \la \jo 0,
    \label{1.2.2} \end{equation}
which is a lax monoidal functor: it is subadditive: $ (\la + \mu)^+ \le \la^+ + \mu^+$, and 
preserves joins (i.e.\ limits).

\Ndt (c) Families of metrics can be viewed as metrics with values in a product of the previous 
structures, like $\Bw^p \ti (\Bw^+)^q$.

\subsection{Multiplicative weights}\label{1.3}
Occasionally, we shall also use the same category $ ([0, \infty], \ge) $ equipped with the monoidal 
structure $\Bw^\bu$, defined by the {\em multiplicative} tensor product $\la\mu$.

	Here we (must) take $ \la\infty = \infty $ for all $ \la$, including $ 0 $ (the initial object 
$\infty $ must be preserved).

	The exponential and logarithm functions give an isomorphism (of monoidal closed 
categories)
    \begin{equation} \begin{array}{c}
\exp\c \Bw   \rlw    \Bw^\bu \cc \ln,
\\[5pt]
e^{\la + \mu} =  e^\la e^\mu,   \qq   e^{\nu - \mu}  =  e^\nu / e^\mu.
    \label{1.3.1} \end{array} \end{equation}

	Accordingly, the internal hom $\, \nu / \mu \,$ of $ \Bw^\bu $ (an extended quotient) gives 
$ \infty/\infty = 0 = 0/0$, applying \eqref{1.1.3}.

\begin{small}

\skp	Again, $ \Bw^\bu $ is *-autonomous, with dualising object the multiplicative identity 1 
and involution $ \la \mapsto \hom^\bu(\la, 1) = 1/\la$.

\end{small}

\section{Real and positive metrics}\label{s2}
	Lawvere's generalised metric spaces \cite{L1} have been used in \cite{G3, G4} as a first 
setting for weighted algebraic topology. Here we deal with the real case.

	Both definitions derive from the formal theory of enriched categories \cite{EK, Ky}, which has 
a simple expression, by inequalities, when the base of enrichment is a category of order-type, 
like $ \Bw^+ $ or $ \Bw$. Distances in a general coquantale were introduced by R. Flagg \cite{Fl}, 
in 1997.

\subsection{Real valued metrics}\label{2.1}
(a) A {\em real metric space} $X$, or $\rho$-{\em metric space}, or $\rho$-{\em space}, 
will be a set $X$ equipped with a $\rho$-{\em metric}, that is, a mapping 
$\rho\c X \ti X \to \Bw$ satisfying the axioms
    \begin{equation}
\rho(x, y) + \rho(y, z)  \ge  \rho(x, z),   \q\;\;   0  \ge  \rho(x, x)   \q   (x, y, z \in X).
    \label{2.1.1} \end{equation}

	We shall generally write this structure $ (X, \rho) $ as $ X$, and the underlying set as 
$ |X|$; the mapping $ \rho $ is written as $\rho_X$ when useful.

	The second condition above can be rewritten as: $\rho(x, x) = 0 $ {\em or} $ \Oinfty$, 
since the properties $0 \ge \rho(x, x)$ and $\rho(x, x) + \rho(x, x) \ge \rho(x, x)$ are only satisfied 
by these extended values.

\skp
\begin{small}

	The value $ \Oinfty $ is rather exceptional in $\rho$-spaces, as shown in \ref{3.4}(c), but 
excluding it would make bigger problems: a chaotic $\rho$-metric only takes the value $\Oinfty$, 
see \ref{3.4}(b).

\end{small}

\Ndt (b) Formally, a $\rho$-space amounts to a small category $ X $ enriched over the symmetric 
monoidal closed category $ \Bw $ considered above, with a hom-value $ X(x, y) $ in 
$ [\Oinfty, \infty]$, written as $ \rho(x, y)$. This extended value will be viewed as {\em the cost of 
the transition from $ x $ to $ y$}, while $ \ga(x, y) = - \rho(x, y) $ represents the corresponding 
{\em gain}, or {\em profit}. (The function $ \ga $ is a `real antimetric', satisfying a reverse 
triangular inequality.)

\Ndt (c) The standard examples are related to the {\em standard $\rho$-line} $\rho\bbR$, with its 
cartesian and tensor powers
    \begin{equation}
\rho\bbR  =  (\bbR, \rho),   \qq   \rho(x, y)  =  y - x.
    \label{2.1.2} \end{equation}

\Ndt (d) $\rho\Mtr $ will denote the category of $\rho$-spaces, with {\em cost-reducing maps} 
$f\c X \to Y $
    \begin{equation}
\rho_X(x, x')  \ge  \rho_Y(f(x), f(x')),   \q   \text{for all } x, x' \in X.
    \label{2.1.3} \end{equation}

\skp
\begin{small}

	Formally, $f$ is a \Bw-enriched functor $X \to Y$, with hom-map $ X(x, x') \to Y(f(x), f(x'))$.

\end{small}

\skp	These mappings will also be called {\em 1-Lipschitz maps} of $\rho$-metric spaces, or 
$\rho$-{\em maps}. Isomorphisms in this category are isometric, and will be called 
{\em isometric isomorphisms} or {\em 1-Lipschitz isomorphisms}. More general `Lipschitz maps' 
will be introduced in \ref{2.5}.

	The {\em reversor} endofunctor $ R\c \rho\Mtr \to \rho\Mtr $ is based on the {\em opposite} 
$\rho$-metric space $ R(X) = X\op$, with the opposite $\rho$-metric, $ \rho\op(x, y) = \rho(y, x)$. 

\subsection{An overview of the main properties}\label{2.2}
The following properties of a $\rho$-metric space $ X = (X, \rho) $ will be studied below.

\Ndt (a) {\em Positivity and finiteness}. A {\em positive} $\rho$-metric, or $\de$-metric, takes 
values in the coquantale $\Bw^+$; this case is reviewed in \ref{2.4}. A {\em finite} $\rho$-metric 
takes values in $\bbR$.

\Ndt (b) {\em Flatness}. We say that a point $ x \in X $ is {\em regular} (resp.\ {\em flat}) if 
$ \rho(x, x) = 0 $ (resp.\ $\Oinfty$). The space itself is {\em regular} (resp.\  {\em flat}) if all its 
points are. Every $\rho$-space $ X $ has a {\em flat part} $ X_\Oinfty$, formed by its flat points. 
A map of $\rho$-spaces can be restricted to the flat parts.

	All important examples here are regular. Flat metrics include the chaotic ones. It is easy 
to see that a $\rho$-space $X$ is flat if and only if $\rho(x, x') = \pm\infty$, for all $x, x' \in X$.

\skp
\begin{small}

	Indeed, if $X$ is flat, $\rho(x, x') \le \Oinfty + \rho(x, x')$, whence $\rho(x, x') = \pm\infty$. 
The converse is also obvious, as the only possible values of $\rho(x, x)$ are 0 and $\Oinfty$.

\end{small}

\Ndt (c) {\em Linearity and invariance}. There are two important classes of real metrics. 
The standard $\rho$-line $ \rho\bbR$, introduced in \ref{2.1}(c), belongs to both.

	First, a {\em linear} $\rho$-metric will be a {\em finite} mapping $\rho\c X \ti X \to \bbR$ 
that satisfies the `additive property' $ \rho(x, y) + \rho(y, z) = \rho(x, z)$; in particular, 
$\rho(x, x) = 0$ for all $x$. We shall see in Section \ref{3.2} that a linear $\rho$-metric has the form 
$\, \De\Phi(x, y) = \Phi(y) - \Phi(x)$, for an arbitrary `potential' function $ \Phi\c X \to \bbR$. Many 
examples related to diverse sciences, can be obtained in this way, using for instance the 
gravitational potential determined by a distribution of maps or the electric potential determined 
by a distribution of charges, or the variation of any variable $ \Phi $ (see \ref{0.2} and \ref{6.6}).

\begin{small}

	Any chaotic metric is additive, in this sense; without the finiteness condition this is 
of less interest.

\end{small}

\skp	Secondly, a $\rho$-metric defined on an abelian group can be {\em invariant under 
translations}, or -- equivalently -- defined by a `real norm' (see \ref{3.3}). Most of the metrics 
we shall consider on the line and its cartesian or tensor products will be of this kind.

\Ndt (d) {\em Symmetry}. A {\em symmetric} $\rho$-space has $ \rho = \rho\op$. They form a 
full subcategory $ !\rho\Mtr \sub \rho\Mtr$, which is reflective and coreflective, see \ref{4.1} 
and \ref{4.2}.

\skp
\begin{small}

	Proposition 3.6 of \cite{G5} shows that a symmetric $\rho$-space is the sum of a positive 
part and a flat part (which is empty, in the important examples).

\end{small}

\skp	More generally, according to a general definition in Directed Algebraic Topology 
(\cite{G4}, 1.2.1), the $\rho$-metric space $ X$ is {\em reversive} if it is isometrically isomorphic 
to its opposite. This is the case of $ \rho\bbR $ and its cartesian or tensor products.

\subsection{Comparisons, limits and colimits}\label{2.3}
(a) The notation $ X \ge X' $ will mean that these $\rho$-spaces have the same underlying set 
and $ \rho_X \ge \rho_{X'}$, or equivalently that the identity of the underlying set is a 
$\rho$-map $ X \to X'$, called a {\em comparison}. We also say that $ \rho_X $ is {\em finer} 
than $ \rho_{X'}$.

	The $\rho$-metrics on a given set form a complete lattice: given a family of them $(\rho_i)$, 
their pointwise join $ \rho(x, y) = \Jo \rho_i(x, y) $ is easily seen to be a $\rho$-metric.

	The category $ \rho\Mtr $ has all limits (resp.\ colimits), computed as in $ \Set $ and 
equipped with the least (resp.\ greatest) $\rho$-metric that makes the structural mappings 
into $\rho$-maps. We write the cartesian power $(\rho\bbR)^n$ as $ \rho\bbR^n$.

\Ndt (b) A product $\, \Pro_i X_i \,$ has the $l_\infty$-type $\rho$-metric
    \begin{equation}
\rho(\Ux, \Uy)  =  \Jo \rho_i(x_i, y_i)   \qq   (\Ux = (x_i), \;  \Uy = (y_i)),
    \label{2.3.1} \end{equation}
and the terminal object is the $\rho$-{\em singleton} $ \rho\top$, equipped with the least 
$\rho$-metric, the chaotic one
    \begin{equation}
\rho\top  =  \sing,   \q   \rho(*, *)  =  \Oinfty   \qq   \text{(terminal object).}
    \label{2.3.2} \end{equation}

	Equalisers have the restricted $\rho$-metric.
	
\Ndt (c) A sum $\, \Sum_i \, X_i$ has the obvious $\rho$-metric (and the initial object is $\es$)
    \begin{equation}
\rho((x, i), (y, i))  =  \rho_i(x, y),   \q   \rho((x, i), (y, j))  =  \infty   \q   (i \neq j).
    \label{2.3.3} \end{equation}

	Coequalisers have the $\rho$-metric induced on the coequaliser in $ \Set$, a quotient set.

	Explicitly, the quotient $ X/R $ of a $\rho$-metric space has the greatest $\rho$-metric 
making the canonical projection a $\rho$-map
    \begin{equation} \begin{array}{c}
\rho(\xi, \eta)  =  \inf_\Ux \; \Sum_j \, \rho(x_{2j-1}, x_{2j}),
\\[6pt]
\Ux  =  (x_1,..., x_{2n})  \in  X^{2n},   \q\;\;   x_1 \in \xi, \;\;  x_{2n} \in \eta,
\\[4pt]
x_{2j} R \, x_{2j+1} \; \text{ for } \; j = 1, ..., n-1.
    \label{2.3.4} \end{array} \end{equation}
%
%

\subsection{Positive directed metrics}\label{2.4}
We recall the basic points of the generalised metrics introduced in \cite{L1}. They were used as 
a first setting for Weighted Algebraic Topology in \cite{G3, G4}.

\Ndt (a) Following the terminology that we used in \cite{G3, G4}, a $\de$-{\em metric space} 
$ X$, or $\de$-{\em space}, or {\em positive} $\rho$-{\em space}, is a set $ X $ equipped with 
a $\de$-{\em metric} $ \de\c X \ti X \to [0, \infty]$, also called a {\em distance}, satisfying the 
axioms
    \begin{equation}
\de(x, y) + \de(y, z)  \ge  \de(x, z),   \qq    0 \ge \de(x, x).
    \label{2.4.1} \end{equation}

	Again, $ X $ is a category enriched over the symmetric monoidal closed category 
$ \Bw^+$, with hom-value $ X(x, y) = \de(x, y) \in [0, \infty]$. We already remarked that  
$\de(x, x) = 0$, for every $x \in X$.

\Ndt (b) $ \de\Mtr $ denotes the full subcategory of $ \rho\Mtr $ of $\de$-metric spaces. Their 
(cost-reducing) maps are also called {\em 1-Lipschitz maps} of $\de$-spaces, or 
$\de$-{\em maps}.

	The $\de$-metrics on a given set form a complete lattice; all limits and colimits exist in 
$ \de\Mtr$. The terminal object is now the $\de$-{\em singleton}, equipped with the least 
$\de$-metric
    \begin{equation}
\de\top  =  \sing,   \qq   \de(*, *)  =  0.
    \label{2.4.2} \end{equation}

	All non-empty products, all equalisers and all colimits are calculated as we have seen in 
$ \rho\Mtr$.

\Ndt (c) The embedding $ \de\Mtr \to \rho\Mtr $ has no left adjoint, as it does not preserve the 
terminal. The coreflector comes from the functor $ ( \; )^+\c \Bw \to \Bw^+ $ of \eqref{1.2.2}, 
which truncates at $ 0 $ the negative values of a $\rho$-metric (yielding the smallest 
$\de$-metric greater than the original $\rho$-metric)
    \begin{equation} \begin{array}{ccc}
( \;\; )^+\c \rho\Mtr \to \de\Mtr    &&    \qqq\q  \text{(coreflector),}
\\[5pt]
Y  \mapsto  Y^+,    &&    \rho^+(y, y')  =  (\rho(y, y'))^+  =  \rho(y, y') \jo 0,
\\[5pt]
Y^+  \iso  Y \ti \de\top,    &&    y  \mapsto  (y, *).
    \label{2.4.3} \end{array} \end{equation}

\begin{small}
\skp

	The absolute value does not give a functor from $ \rho\Mtr $ to $ \de\Mtr$: for instance, it 
cannot work on the $\rho$-singleton, nor on any $\rho$-space having some flat point. But it 
does act on linear $\rho$-metrics, see \ref{5.2}(c).

\end{small}

\Ndt (d) Again, the reversor endofunctor $ R\c \de\Mtr \to \de\Mtr $ is based on the {\em opposite} 
$\de$-metric space $ R(X) = X\op$, with the opposite $\de$-metric, $ \de\op(x, y) = \de(y, x)$.

	A {\em symmetric} $\de$-metric $ \de = \de\op $ is the same as an {\em \'ecart} in Bourbaki 
\cite{Bk}. We shall denote by $ !\de\Mtr $ the full subcategory of $ \de\Mtr $ of {\em symmetric} 
$\de$-{\em metric spaces}. (This category was written $ \Mtr $ in \cite{G3, G4}, but here we want 
to be consistent with the parallel case $ !\rho\Mtr$.)

\Ndt (e) {\em Structures on the line and the interval}. The basic example, as in \cite{G3, G4}, is 
the {\em standard }$\de$-{\em line} $ \de\bbR$, with $\de$-metric
    \begin{equation}
\de(x, y)  =  y - x,  \IF  x \le y,   \qq   \de(x, y)  =  \infty,  \text{ otherwise,}
    \label{2.4.4} \end{equation}
where $ y < x $ gives an infinite distance, with the aim of directing paths by precluding 
backward moves. Its cartesian powers are written $ \de\bbR^n$.

	The sub-$\de$-space $ \de\bbI \sub \de\bbR $ (that is, the interval $ [0, 1] $ with the 
restricted $\de$-metric) will play also here, as in \cite{G3, G4}. the role of the standard interval, 
on which paths and homotopies are parametrised: see \ref{5.1}.

	Of course we also use the {\em euclidean line} $\bbR$, with the usual euclidean metric $d$.
	
	A different structure, used in \cite{L1}, is also of interest, the $\de_0$-{\em line} 
$\, \de_0\bbR$
    \begin{equation}
\de_0\bbR  =  (\rho\bbR)^+  =  (\bbR, \de_0),   \q   \de_0(x, y)  =  (y - x) \jo 0,
    \label{2.4.5} \end{equation}
whose cartesian powers are written $ \de_0\bbR^n$.

\subsection{Lipschitz maps}\label{2.5}
(a) Let $\la \in \bbR^+ = [0, \infty[$. We shall write $\la \ge 0$, leaving understood that $\la < \infty$. 

	If $\la > 0$, multiplication by $\la$ in $[\Oinfty, \infty]$ preserves and reflects the order; we  
also assume that $ 0\Oinfty = \Oinfty $ and $ 0\infty = \infty$ (cf.\ \ref{1.3}).

	We say that a mapping $ f\c X \to Y $ between $\rho$-spaces is $\la$-{\em Lipschitz}, or a 
$\la$-{\em map}, if
    \begin{equation}
\la\rho_X(x, x')  \ge  \rho_Y(f(x), f(x'))   \qq   (x, x' \in X),
    \label{2.5.1} \end{equation}
and we say that $ f $ is a {\em Lipschitz map}, or $\rho_\infty$-{\em map}, if it is has some 
{\em Lipschitz constant} $ \la \ge 0$, in the previous sense. Composing a $\la$-map with a 
consecutive $\la'$-map we get a $\la\la'$-map.

	If $ X $ is positive (a $\de$-space) one can replace $ \la $ with any $ \la' \ge \la$, as 
examined in (b). Generally, one cannot: see \ref{2.6}. (Lipschitz maps in particular cases are 
characterised in \ref{3.2}, \ref{4.5}, \ref{4.6}.)

	The category $ \rho\Mtr $ is a wide subcategory of the category $ \rho_\infty\Mtr $ of 
$\rho$-spaces and Lipschitz maps. A {\em Lipschitz isomorphism} will be an isomorphism of 
the latter. The symmetric case is denoted as $ !\rho_\infty\Mtr$.

	Restricting to $\de$-metric spaces, we have the full subcategory $ \de_\infty\Mtr $ of 
$\de$-metric spaces and their Lipschitz maps, also called $\de_\infty$-{\em maps}, and its 
full subcategory $ !\de_\infty\Mtr $ for the symmetric case.

\Ndt (b) {\em If $ X $ is positive} (a $\de$-space) and $ Y $ an arbitrary $\rho$-space, one can 
define, for every mapping $ f\c X \to Y $ (of sets) a {\em Lipschitz weight}, or {\em seminorm} 
$||f|| \in [0, \infty]$
    \begin{equation}
||f||  =  \inf\{\la \in [0, \infty[  \sep  \forall\;  x, x' \in X,  \; \rho_Y(fx, fx') \le \la\de_X(x, x')\},
    \label{2.5.2} \end{equation}
and $ f $ is Lipschitz if and only if $||f|| < \infty$. In this case, $ ||f|| $ is the least Lipschitz 
constant for $ f$, and every $ \la \ge ||f|| $ is also a Lipschitz constant.

	The categories $ \de_\infty\Mtr $ and $ !\de_\infty\Mtr $ are {\em multiplicatively weighted} 
(or seminormed) by the Lipschitz weight. By this we mean that every morphism $ f $ is assigned 
a weight $ ||f|| \in [0, \infty[ $ so that
    \begin{equation}
||gf|| \le ||f|| \; ||g||,   \qq   || \id X || \le 1.
    \label{2.5.3} \end{equation}

	These two axioms imply that the weight of an identity can only be $ 1 $ or $ 0$. The wide 
subcategories $ \de\Mtr $ and $ !\de\Mtr $ are characterised by the condition $ ||f|| \le 1 $ and 
inherit the restricted weight.
	
\Ndt (c) The category $ \rho_\infty\Mtr $ (resp.\ $ !\rho_\infty\Mtr$, $ \de_\infty\Mtr$, 
$ !\de_\infty\Mtr$) is {\em finitely} complete and cocomplete, and the inclusion of $ \rho\Mtr $ 
(resp.\ $ !\rho\Mtr$, $\de\Mtr$, $!\de\Mtr$) preserves finite limits and colimits. But the `metric' of 
a (co)limit in the larger category is only determined up to Lipschitz isomorphism.

	If $ X $ is a $\rho$-metric space and $ \la \in [0, \infty[$, we shall write $ \la X $ the same 
set equipped with the $\rho$-metric $ \la\rho_X$. Thus, a $\la$-Lipschitz map is the same as a 
$\rho$-map $ \la X \to Y$. This transformation preserves symmetric structures and positive 
structures.

\subsection{Remarks}\label{2.6}
The extension to real valued metrics brings about aspects which strongly diverge from the positive 
case. We begin with some elementary observations. In particular, they show that the structure of 
$\rho\bbR$ highly restrains the Lipschitz maps {\em defined} on this space; in fact, we are 
more interested in maps {\em from $\de$-spaces to} $ \rho\bbR $ (as in \ref{4.6}).

\Ndt (a) A $\rho_\infty$-map $f\c \rho\bbR \to \rho\bbR$ is simply an affine map $f(x) = \la x + y_0$ 
(with $ \la  \ge 0$), and admits a unique Lipschitz constant, namely $ \la$. (This will be extended 
to linear $\rho$-spaces, in \ref{3.2}(d).)

	In fact, interchanging $ x $ and $ x'$, the condition $ \la(x' - x) \ge fx' - fx $ is equivalent to 
$ \la(x' - x) = fx' - fx$, that is, $ f(x) = \la x + f(0)$.

\Ndt (b) Any $\la$-Lipschitz maps $ f\c \rho\bbR \to \de_0\bbR $ is constant, and $ \la $ is 
necessarily 0 (see \ref{3.2}(e)).

\begin{small}
\skp

	Letting $ t > t' $ in $ \bbR$, we have $ \la(t' - t) \ge \de_X(ft, ft') \ge 0$, which gives 
$ \de_X(ft, ft') = 0 $ for $ \la = 0 $ and is impossible for $ \la > 0$.

\end{small}

\Ndt (c) {\em Points as maps}. A point of a $\rho$-space $ X $ can be identified with a 
$\rho$-map $ \de\top \to X $ {\em defined on the} $\de$-{\em singleton}. Formally, the forgetful 
functors $ \rho_\infty\Mtr \to \Set $ and $ \rho\Mtr \to \Set $ are represented by the 
$\de$-singleton (the identity of the tensor product, see \ref{6.1}). If $ X $ is positive every 
constant mapping $f\c X \to Y$ is a $\rho$-map, since it can be factorised as $X \to \de\top \to Y$.

	On the other hand, a $\rho_\infty$-map $ \rho\top \to Y $ can only reach a flat point of 
$ Y$. (The $\rho$-singleton represents the `flat-part functor' mentioned in \ref{2.2}(b).)

\subsection{Future and past topologies}\label{2.7}
A $\rho$-space $ X = (X, \rho) $ can be equipped with several topologies, which coincide in 
the symmetric case.

\Ndt (a) The {\em future topology} $\, \tp^+(X) $ is defined by a canonical local base at each point 
$ x_0 \in X$, the {\em future} open $\ep$-balls 
    \begin{equation}
B^+(x_0, \ep)  =  \{x \in X  \sep  \rho(x_0, x) < \ep\}     \q\;\;     (\ep > 0  \inn  \bbR),
    \label{2.7.1} \end{equation}
also written as $ B^+_X(x_0, \ep)$. (Of course, this local base can be made countable.) A finer 
$\rho$-metric $ \rho' \ge \rho $ (in \ref{2.3}(a)) gives a finer future topology, with more open sets.

	Dually, the {\em past topology} $\, \tp^-(X) = \tp^+(X\op)$ is defined by the system of 
{\em past} open $\ep$-balls
    \begin{equation}
B-(x_0, \ep)  =  \{x \in X  \sep  \rho(x, x_0) < \ep\}     \q\;\;     (\ep > 0  \inn  \bbR).
    \label{2.7.2} \end{equation}

	However, we shall more often use the `reflective' and `coreflective' symmetric topologies, 
defined in \ref{4.1} and \ref{4.2}, keeping a trace of the asymmetry of $ \rho $ by a preorder (see  \ref{4.3}).

\Ndt (b) For a symmetric $\rho$-space $ X $ the {\em associated topological space} 
$\, \tp(X) = \tp^+(X) = \tp^-(X) $ will generally be written $ X$.

	We have defined three functors
    \begin{equation} \begin{array}{c}
\tp^+\c \rho_\infty\Mtr \to \Top,   \q   \tp^-\c \rho_\infty\Mtr \to \Top,
\\[5pt]
\tp\c !\rho_\infty\Mtr \to \Top,
    \label{2.7.3} \end{array} \end{equation}
and will use the same notation for their restriction to the positive case.

\section{Basic examples and main properties}\label{s3}
	After examining $\rho$-structures on the line and recalling $\de$-structures on the 
spheres, we deal with the main properties we are interested in: linear metrics, defined by a 
potential function (in \ref{3.2}), and invariant metrics, defined by a real norm (in \ref{3.3}).

	Several more technical properties of flat and symmetric $\rho$-spaces can be found in 
the preprint version of this article, Sections 3.5-3.6 \cite{G5}.

\subsection{Metrics on the line and the spheres}\label{3.1}
(a) We shall generally test our notions on three basic structures, from which many others will be 
deduced later: the $\de$-line $ \de\bbR $ (in \eqref{2.4.4}), the $\rho$-line $ \rho\bbR $ (in 
\ref{2.1}(c)), and the Lawvere directed line $ \de_0\bbR $ (in \eqref{2.4.5}). Their future topology 
is determined by the following canonical local base at any point $ x_0$:
    \begin{equation*} \begin{array}{lll}
\text{- for }  \de\bbR\c &   B^+(x_0, \ep)  =  [x_0, x_0 + \ep[,
\\[5pt]
\text{- for }  \rho\bbR  \and  \de_0\bbR\c  \;\;\;  &  
B^+(x_0, \ep) = \; ]\Oinfty, x_0 + \ep[  \q   &  (\ep > 0).
    \label{} \end{array} \end{equation*}

	This begins to make clear why we are using $ \de\bbR $ as the standard (metric) 
{\em directed} line: the `cone of the future' at $ x_0 $ is appropriate to our goals, while that of 
$ \de_0\bbR $ is not. The $\rho$-line $ \rho\bbR $ suits different goals, and its paths will be 
allowed to undergo losses and gains (as in the example of \ref{0.2}).

	We also note that the previous metrics and topologies are invariant under translation.

\Ndt (b) The orbit space $ (\de\bbR)/\bbZ $ (with respect to the action of the group $\bbZ$ on 
the line, by integral translations) gives the standard $\de$-circle $\de\bbS^1$, where $\de(x, y)$ 
is the measure of the counterclockwise arc from $x$ to $ y$, with respect to the length of the 
circle \cite{G3, G4}. The higher $\de$-spheres also have an interesting structure, a quotient of 
the $\de$-cube $ \de\bbI^n$.

\Ndt (c) On the other hand, the orbit space $ (\rho\bbR)/\bbZ $ has the chaotic $\rho$-metric, 
always $ \Oinfty$, because two points $ [x], [y] $ can be lifted to $ \rho\bbR $ with arbitrarily 
small $ \rho(x, y)$.

	Relative metric spaces lack a `standard $\rho$-circle'. As already said in the Introduction, 
we shall introduce a more general setting in Part III.

\begin{small}
\skp

	As an {\em ad hoc} solution, one might introduce a notion of `local $\rho$-metric', much in 
the same way as the chaotic preorder induced by the ordered line on the circle can be replaced 
by a local preorder (see \cite{G4}, 1.9.3). Yet, this complicated framework would be ineffective for 
other drawbacks of real metrics.

\end{small}

\subsection{Linear real metrics and potential functions}\label{3.2}
Let $ (X, \rho) $ be a $\rho$-metric space.

\Ndt (a) The $\rho$-metric is linear (see \ref{2.2}) if and only if it is finite and 
$ \rho(x, y) = - \rho(y, x)$, for all pairs $x, y$. Moreover, $ \rho(x, x) = 0$, for all $ x$.

\begin{small}
\skp

	Indeed, if $ \rho $ is linear (which includes being finite, by definition), $ \rho(x, y) + \rho(y, x) 
= \rho(x, x) = 0$. Conversely, assuming that $ \rho(x, y) = - \rho(y, x)$, the triangular inequality 
implies the reverse one.

	Let us note that the condition $ \rho(x, y) = - \rho(y, x) $ does not imply that our $\rho$-metric 
is finite: for instance, on the set $ \{0, 1\} $ one can take $ \rho(0, 0) = 0 = \rho(1, 1) $ and 
$ \rho(0, 1) = \infty = - \rho(1, 0)$.

\end{small}

\Ndt (b) The $\rho$-metric is linear if and only if there is a function $ \Phi\c X \to \bbR $ 
(a {\em potential} of $ \rho$) such that $ \rho $ is its {\em variation} $\De\Phi$
    \begin{equation}
\rho(x, y)  =  \De\Phi(x, y)  =  \Phi(y) - \Phi(x).
    \label{3.2.1} \end{equation}

\begin{small}

	This expression plainly gives a linear $\rho$-metric. Conversely, if $\rho$ is linear, we 
define $ \Phi(x) = \rho(x_0, x)$, for a fixed point $ x_0 \in X$, and verify \eqref{3.2.1}, using (a): 
$\Phi(y) - \Phi(x)  =  \rho(x_0, y) - \rho(x_0, x)  =  \rho(x_0, y) + \rho(x, x_0)  =  \rho(x, y).$

\end{small}

\skp	Concretely, $ \Phi $ could be the measure of any `variable' on the `space' $ X$. Of course, if 
$\Phi$ is a $C^1$-function on a differentiable manifold, it is the potential of its differential 1-form. 
One can view the variation $\rho = \De\Phi$ as a `finite differential' of $\Phi$.

	The potential $ \Phi $ is determined up to an additive constant (for $ X\neq \es$, of course): 
fixing an arbitrary point $ x_0 \in X$, we can always take $ \Phi(x_0) = 0 $ and
    \begin{equation}
\Phi(x)  =  \rho(x_0, x).
    \label{3.2.2} \end{equation}

\Ndt (c) Obvious examples can be deduced from the gravitational or electric fields, as in \ref{0.2} 
and \ref{6.6}. A non-conservative field would produce a more general structure, studied in Part III.

	A product of linear $\rho$-spaces is not linear, generally: for instance $ \rho\bbR^2 $ is not. 
A linear $\de$-metric is discrete: $ \de(x, y) = 0$, for all $ x, y$.

\Ndt (d) If $ X $ and $ Y $ are linear $\rho$-spaces, any $\rho$-map $ f\c X \to Y $ is an isometry. 
Any $\la$-Lipschitz map is an isometry $\la X \to Y$, and $\la$ is uniquely determined, 
unless $ X $ has the chaotic $\de$-metric, everywhere 0.

\begin{small}

\skp	In fact, from $ \la\rho_X(x, x') \ge \rho_Y(fx, fx') $ and (a) we deduce that $ \la\rho_X(x, x')$ 
$= \rho_Y(fx, fx')$, so that $ f $ is an isometry $ \la X \to Y$. Since linear $\rho$-metrics are finite, 
$ \la $ is uniquely determined unless $ \rho_X$ is everywhere $ 0 $ (which includes the case 
$ X = \es$).

\end{small}

\Ndt (e) If $ f\c X \to Y $ is a $\la$-Lipschitz map from a linear $\rho$-space to a $\de$-space, then

\Ndt - $ X $ has the chaotic $\de$-metric or $ \la = 0 $ (possibly both),

\ndt - $ f(X) $ has the chaotic $\de$-metric.

\skp	When these conditions are satisfied, $ f $ is an arbitrary mapping.

\begin{small}

\skp	In fact, applying (a), the condition $ \la\rho_X(x, x') \ge \de_Y(fx, fx') \ge 0 $ is equivalent to 
$ \la\rho_X(x, x') = 0 = \de_Y(fx, fx')$.

\end{small}

\subsection{Invariance and real norms}\label{3.3}
On a set $ X $ acted on by a group $ G $ we can consider a $\rho$-metric invariant up to the 
action of the group. More simply, we suppose that $ X $ is an abelian group in additive notation, 
acting on itself.

\Ndt (a) A $\rho$-metric $ \rho $ on $ X $ is said to be {\em invariant up to translations} if
    \begin{equation}
\rho(x + z, y + z)  =  \rho(x, y)   \qqq   (x, y, z \in X).
    \label{3.3.1} \end{equation}

	Here $ \rho(x, x) = \rho(0, 0)$, for all $ x$: the metric is {\em either} regular {\em or} flat, 
according to the definitions of \ref{2.2}(b) -- also because $X$ is non-empty. In the second (less 
interesting) case, we have already seen that $\rho(x, y) = \pm \infty$, for all $ x, y \in X$.

\Ndt (b) One can equivalently assign a {\em real norm}, or $\rho$-{\em norm}, on $ X$, 
which means a mapping $\Umu\c X \to \Bw$ satisfying the axioms
    \begin{equation}
\Umu(x) + \Umu(y)  \ge  \Umu(x + y),   \q   \Umu(0)  \le  0   \q\;\;   (x, y \in X),
    \label{3.3.2} \end{equation}
where $ \Umu(0) $ can only be $ 0 $ or $ \Oinfty$.

	The correspondence between metric and norm is obviously expressed by the following 
relations:
    \begin{equation}
\Umu(x)  =  \rho(0, x),   \qq   \rho(x, y)  =  \Umu(y - x).
    \label{3.3.3} \end{equation}

\begin{small}

\skp	In fact, if $ \rho\c X \ti X \to \Bw $ is an invariant $\rho$-metric, the associated mapping 
$ \Umu\c X \to \Bw $ satisfies \eqref{3.3.2}
    \begin{equation*} \begin{array}{l}
\Umu(x) + \Umu(y) = \rho(0, x) + \rho(0, y) = \rho(0, x) + \rho(x, x + y)
\\[3pt]
\qq\;\;\;\;\;   \ge \; \rho(0, x + y) = \Umu(x + y).
    \label{} \end{array} \end{equation*}

	Conversely, assuming that $ \Umu\c X \to \Bw $ satisfies \eqref{3.3.2}, the associated 
mapping $ \rho\c X \ti X \to \Bw $ is obviously invariant up to translations and
    \begin{equation*}
\rho(x, y) + \rho(y, z) = \Umu(y - x) + \Umu(z - y) \,  \ge  \, \Umu(z - x) =  \rho(x, z).
    \label{} \end{equation*}

\end{small}

\Ndt(c) A product of abelian groups with $\rho$-metrics invariant up to translations is of the 
same kind. For instance all $ \de\bbR^n$, $\rho\bbR^n$, $\de_0\bbR^n $ have a 
$\rho$-metric invariant up to translations.

\Ndt (d) A $\rho$-metric $ \rho $ on $ X $ which is invariant up to translations {\em and} linear 
(as in $ \rho\bbR$) amounts to a homomorphism of abelian groups
    \begin{equation}
\Phi\c  X \to \bbR,   \qq   \Phi(x + y)  =  \Phi(x) + \Phi(y),
    \label{3.3.4} \end{equation}
which is at the same time its potential (null at $ 0_X$) and its real norm
    \begin{equation}
\rho(x, y)  =  \Phi(y) - \Phi(x)  =  \Phi(y - x).
    \label{3.3.5} \end{equation}

\begin{small}

\skp	In fact, if $ \rho $ is invariant up to translations and linear, $ \Umu(x) = \rho(0, x) = \Phi(x)$, 
and the relation \eqref{3.3.5} says that $ \Phi $ is a homomorphism. Conversely, \eqref{3.3.5} 
trivially implies that $ \Phi $ is a $\rho$-norm for $ \rho$.

\end{small}

\subsection{Discrete and chaotic metrics}\label{3.4}
(a) Every set $ X $ can be equipped with the {\em discrete} $\rho$-{\em metric} $ \rho_D $ 
(the greatest one), which is positive and symmetric
    \begin{equation}
\rho_D(x, x)  =  0,   \q   \rho_D(x, y)  =  \infty   \qq   (\for x \neq y),
    \label{3.4.1} \end{equation}
giving the {\em discrete} $\rho$-{\em metric space} $ \rho_DX$, and the left adjoint to the 
forgetful functor $ \rho\Mtr \to \Set. $ In $ \de\Mtr $ these issues will be written $ \de_D $ 
and $ \de_DX$. (Of course the associated topology is discrete.)

	A discrete $\rho$-metric space is a sum of $\de$-singletons.

\Ndt (b) Every set $ X $ has a {\em chaotic} $\rho$-{\em metric} and a {\em chaotic} 
$\de$-{\em metric} (the least ones), which are symmetric
    \begin{equation}
\rho_C(x, y)  =  \Oinfty,   \q   \de_C(x, y)  =  0   \qq   (x, y \in X).
    \label{3.4.2} \end{equation}

	The former gives the {\em chaotic} $\rho$-{\em metric space} $ \rho_CX$, and the right 
adjoint to the forgetful functor $ \rho\Mtr \to \Set$. The latter gives the {\em chaotic} 
$\de$-{\em metric space} $ \de_CX$, and the right adjoint to $ \de\Mtr \to \Set$. (In both cases 
the associated topology is chaotic.)

\Ndt (c) We say that a $\rho$-space is {\em affordable} if $ \rho(x, x') < \infty$, for all pairs of 
points. The same term applies to $\de$-spaces.

	An affordable $\rho$-space $ X $ that has some points $ z, z' $ with $\rho(z, z') = \Oinfty$ 
is chaotic: $ \rho(x, y) \le \rho(x, z) + \rho(z, z') + \rho(z', y) = \Oinfty$, for all $ x, y $ (since all 
values are $ < \infty$).

\section{Symmetrisation, topology and preorder}\label{s4}

	We show that the full subcategory $ !\rho\Mtr \sub \rho\Mtr $ of symmetric $\rho$-metric 
spaces (see \ref{2.2}(d)) is reflective and coreflective in $ \rho\Mtr$, producing two associated 
`symmetric' topologies. In the positive case, the reflector was used in \cite{G3, G4} while the 
coreflector was used in \cite{L1}. Similarly, a $\rho$-metric space has a reflective and a 
coreflective preorder, defined in \ref{4.3}.

	Both these aspects will be of use here, in different analyses -- as shown in \ref{4.5}. 
Their relationship to cartesian and tensor products is dealt with in Section 6.

	The `indeterminate' relation of real valued metrics to topology, discussed in the Introduction, 
is also investigated in \ref{4.8}.

\subsection{The reflective symmetric metric and topology}\label{4.1}
(a) The reflector $^\me\c \rho\Mtr \to \; !\rho\Mtr $ is based on the greatest symmetric $\rho$-metric $ 
\hat{\rho} \le \rho $ (already used in \cite{G3, G4} in the positive case, and written $!\de$), which 
will be called here the {\em reflective symmetric} $\rho$-metric of $ \rho $
    \begin{equation} \begin{array}{c}
^\me\c \rho\Mtr \to \; !\rho\Mtr,  \q   (X, \rho)^\me  =  (X, \hat{\rho}),  \q   \hat{\rho}  =  \rho \me \rho\op,
\\[5pt]
\hat{\rho}(x, x')  =  \inf_\Ux (\Sum_j \, (\rho(x_{j-1}, x_j) \me \rho(x_j, x_{j-1}))),
\\[5pt]
\hat{\rho}(x, x')  \le  \rho(x, x') \me \rho(x', x),
    \label{4.1.1} \end{array} \end{equation}
where $ \Ux = (x_0,..., x_n) $ is any finite sequence in $ X $ with $x_0 = x$ and $x_n = x' $ (a 
{\em step-path} from $ x $ to $ x'$). The unit is the comparison $X \to X^\me$, the counit is the 
identity $ Y^\me = Y$, for a symmetric $\rho$-space. (The inclusion functor 
$ !\rho\Mtr \sub \rho\Mtr $ is left understood.)

	Restricting to the positive case, we obtain the reflector $ ^\me\c \de\Mtr \to \; !\de\Mtr$, based 
on the greatest symmetric $\de$-metric $ \de^\me \le \de$. One can replace $ \rho\Mtr $ by 
$ !\rho_\infty\Mtr$, and $ \de\Mtr $ by $ !\de_\infty\Mtr$.

\begin{small}

\skp	Loosely speaking, \eqref{4.1.1} `adjusts' the symmetric function 
$ \rho(x, x') \me \rho(x', x) $ to satisfy the triangular inequality in $ X$. However, {\em if this 
function already satisfies it}, the meet $\rho \me \rho\op$ can be computed pointwise on $X \ti X$
    \begin{equation}
\hat{\rho}(x, x')  =  \rho(x, x') \me \rho(x', x).
    \label{4.1.2} \end{equation}

\end{small}

\Ndt (b) Composing with the functor $ \tp\c !\rho_\infty\Mtr \to \Top $ of \eqref{2.7.3} we have a 
functor $ \tp^\me$, also written $^\me$
    \begin{equation}
^\me = \tp^\me\c \rho_\infty\Mtr \to \Top   \qq   (^\me\c \rho\Mtr \to \Top).
    \label{4.1.3} \end{equation}
which equips a $\rho$-metric space with the {\em reflective (symmetric) topology}, associated to 
$ \hat{\rho}$.

\Ndt(c)  For the basic examples of \ref{3.1}(a)

\LL

\Ndt  - $ (\de\bbR)^\me = \bbR $ has the euclidean metric and topology,

\ndt - $ (\rho\bbR)^\me = \rho_C\bbR $ and $ (\de_0\bbR)^\me = \de_C\bbR $ have a chaotic 
metric and chaotic topology.

\LB
\begin{small}

\skp	In the first case, $\de(x, x') \me \de(x', x) = |x - x'|$ is (already) the euclidean 
metric of the line. The second instance is proved in Lemma \ref{4.7}. In the third 
$\de_0^\me(x, x') = \de_0(x, x') \me \de_0(x', x) = 0$.

\end{small}

	The first case above (and its consequences on cartesian powers, tensor powers and 
substructures) shows that the reflector is appropriate for the goals of \cite{G3, G4}. The 
second shows that we cannot rely on it to get a useful topology on $\rho\bbR$.

\subsection{The coreflective symmetric metric and topology}\label{4.2}
(a) The coreflector $^\jo\c \rho\Mtr \to \; !\rho\Mtr $ is based on the least symmetric $\rho$-metric 
$\check{\rho}\ge \rho$, which will be called the {\em coreflective symmetric} $\rho$-metric of $\rho$
    \begin{equation} \begin{array}{ccc}
^\jo\c \rho\Mtr \to \; !\rho\Mtr,  &\;\;&   (X, \rho)^\jo  =  (X, \check{\rho}),
\\[5pt]
\check{\rho}  =  \rho \jo \rho\op,   &&   \check{\rho}(x, x')  =  \rho(x, x') \jo \rho(x', x).
    \label{4.2.1} \end{array} \end{equation}

	The counit is the comparison $ X^\jo \to X$.

	Let us note that $ \check{\rho} $ is computed pointwise on $ X \ti X $ because $\rho$-metrics 
are closed under joins. Moreover
    \begin{equation}
\hat{\rho}  \le  \rho  \le  \check{\rho},   \q   \hat{\rho}  \le  \rho\op  \le  \check{\rho}.
    \label{4.2.2} \end{equation}

	Restricting to the positive case, we obtain the coreflector $ ^\jo\c \de\Mtr \to \; !\de\Mtr$, based 
on the least symmetric $\de$-metric $ \de^\jo \ge \de$. Again, one can replace $ \rho\Mtr $ by 
$ !\rho_\infty\Mtr$, and $ \de\Mtr $ by $ !\de_\infty\Mtr$.

	These right adjoints preserve products, and all limits.

\Ndt (b) Composing with the functor $ \tp\c !\rho_\infty\Mtr \to \Top $ of \eqref{2.7.3} we have a 
functor $ \tp^\jo$, also written $^\jo$
    \begin{equation}
^\jo = \tp^\jo\c \rho_\infty\Mtr \to \Top   \q   (^\jo\c \rho\Mtr \to \Top),
    \label{4.2.3} \end{equation}
which equips a $\rho$-metric space with the {\em coreflective (symmetric) topology}, associated 
to $ \check{\rho}$. For a $\rho$-space $ X $ we have the following `finer' relations, for metrics 
and topologies:
    \begin{equation}
X^\jo  \ge  \, X,  X\op  \, \ge  X^\jo,   \q\;\;   X^\jo  \ge \,  \tp^+X,  \tp^-X  \, \ge  X^\jo.
    \label{4.2.4} \end{equation}

\Ndt (c) For the basic examples of \ref{3.1}(a)

\LL

\ndt - $ \de\bbR^\jo = \de_D\bbR $ has the discrete $\de$-metric, infinite out of the diagonal, and 
the discrete topology,

\Ndt - $ \rho\bbR^\jo $ and $ \de_0\bbR^\jo $ have the euclidean metric 
$ \check{\rho}(x, x') = \de_0^\jo(x, x') = |x - x'| $ and the euclidean topology.

\LB

\skp	The coreflector is thus of interest for $ \rho\bbR$, and of little use for $ \de\bbR$. We shall 
find convenient ways to use both the reflector and the coreflector (see \ref{4.5}); or none of them, 
as far as topology is concerned (in Part III).

\Ndt (d) The results of this section, up to now, are abridged in the following diagrams:
    \begin{equation} \begin{array}{c} 
    \xymatrix  @C=30pt @R=20pt
{
~!\rho\Mtr~~   \ard[r]    &   ~~ \rho\Mtr~~    \ar@<-4pt>[l]_-{\me}   \ar@<4pt>[l]^-{\jo}   &
~!\rho_\infty\Mtr~~   \ard[r]    &   
~~ \rho_\infty\Mtr~~    \ar@<-4pt>[l]_-{\me}   \ar@<4pt>[l]^-{\jo}
\\ 
~~!\de\Mtr~~   \ard[r]    \ard[u]    &   ~~\de\Mtr~   \ar@<-4pt>[l]_-{\me}   \ar@<4pt>[l]^-{\jo}  \ard[u]  &
~~!\de_\infty\Mtr~~   \ard[r]    \ard[u]    &   
~~\de_\infty\Mtr~   \ar@<-4pt>[l]_-{\me}   \ar@<4pt>[l]^-{\jo}  \ard[u]
}
    \label{4.2.5} \end{array} \end{equation}

	Each of them contains three commutative squares, formed by full inclusions (the 
dashed arrows), their left adjoints (labelled with $^\me$) and the right ones (with $^\jo$).

\subsection{Metrics and preorders}\label{4.3}
The asymmetry of a $\rho$-metric, overcome by symmetrisation, is also important, and can be 
highlighted by an associated preorder relation. Again there are two main ways: the reflective 
one (of interest for $ \de\bbR$, for instance) and the coreflective one (of interest for $ \rho\bbR $ 
and $ \de_0\bbR$).

\Ndt(a) A $\rho$-space $ X $ can be equipped with the {\em reflective preorder} $ x \prec_\infty x'$
    \begin{equation} \begin{array}{c}
\pr_\infty\c \rho_\infty\Mtr \to \pSet,   \q\;\;\;   \pr_\infty(X, \rho)  =  (X, \prec_\infty)
\\[5pt]
x \prec_\infty x'  \; \IF \;  \rho(x, x') < \infty,
    \label{4.3.1} \end{array} \end{equation}
saying that the transition from $ x $ to $ x' $ is affordable.

	This functor preserves {\em finite} products. It also preserves all colimits, being left adjoint 
to the following embedding of preordered sets in (flat) $\rho$-spaces:
    \begin{equation} \begin{array}{c}
\mt_\infty\c \pSet \to \rho_\infty\Mtr,   \q\;\;   \mt_\infty(Y, \prec)  =  (Y, \rho_\infty),
\\[5pt]
\rho_\infty(y, y')  =  \Oinfty   \IF  y \prec y',  \:  \and \;  \infty  \;  \text{ otherwise.}
    \label{4.3.2} \end{array} \end{equation}

\begin{small}

	For a preordered set $ (Y, \prec)$, a mapping $ f\c (X, \rho) \to (Y, \rho_\infty) $ satisfies 
$ \rho_\infty(fx, fx') \le \la\rho(x, x') $ (for some $ \la \ge 0$) if and only if $ \rho(x, x') < \infty $ 
implies $ \rho_\infty(fx, fx') < \infty$, that is, $ x \prec_\infty x' $ implies $ fx \prec fx'$.

\end{small}

\skp	The reflective preorder is of interest for the line $ \de\bbR$, which gets the natural order 
$ x \le x' $ (as in \cite{G3, G4}). It is of no use for $ \rho\bbR $ and $ \de_0\bbR$, which get the 
chaotic preorder.

\Ndt (b) A $\rho$-space $ X $ can also be equipped with the {\em coreflective preorder} 
$ x \prec_0 x'$
    \begin{equation} \begin{array}{c}
\pr_0\c \rho_\infty\Mtr \to \pSet,   \q\;\;   \pr_0(X, \rho)  =  (X, \prec_0),
\\[5pt]
x \prec_0 x'   \IF  \rho(x', x) \le 0,
    \label{4.3.3} \end{array} \end{equation}
saying that the backward transition, from $x'$ to $x$, gives a gain $\ga(x', x) \ge 0$.

	This functor preserves all limits; in fact it is right adjoint to a different embedding
    \begin{equation} \begin{array}{c}
\mt_0\c \pSet \to \de\Mtr \sub  \rho_\infty\Mtr,   \q\;\;   \mt_0(Y, \prec)  =  (Y, \de_0),
\\[5pt]
\de_0(y, y')  =  0   \IF  y' \prec y,   \and  \infty  \;  \text{ otherwise.}
    \label{4.3.4} \end{array} \end{equation}

\begin{small}

	For a preordered set $ (Y, \prec)$, a mapping $ f\c (Y, \de_0) \to (X, \rho) $ satisfies 
$ \rho(fy, fy') \le \la\de_0(y, y') $ (for some $ \la \ge 0$) if and only if $ \de_0(y, y') = 0 $ implies 
$ \rho(fy, fy') = 0$, that is, $ y' \prec y $ implies $ fy' \prec_0 fy$.

\end{small}

\skp	For instance, the structures $ \rho\bbR $ and $ \de_0\bbR = (\rho\bbR)^+ $ get the natural 
order $ x \le x'$, while $ \de\bbR $ gets the discrete one.
	
\Ndt (c) In a $\rho$-space $ X$, the relation $ x \prec_0 x' $ implies $ x' \prec_\infty x $ ({\em 
note the reversion}). In the symmetric case both preorders are equivalence relations.

\subsection{Topologies and preorders on real metric lines}\label{4.4}

The previous results for the basic examples of \ref{3.1}(a) are summarised in this table
    \begin{equation} \begin{array}{cccccc}
(X, \rho)  &  \tp^+X  &  \hat{\rho}(x, x')  &  x \prec_\infty x'  &  \check{\rho}(x, x')  &  
x \prec_0 x'
\\[5pt]
\de\bbR  &  [x_0, x_0 + \ep[  &  |x - x'|  &  x \le x'  &  discrete  &  discrete
\\[5pt]
\rho\bbR, \de_0\bbR  &  [\Oinfty, x_0 + \ep[  &  chaotic  &  chaotic  &  |x - x'|  &  x \le x' 
    \label{4.4.1} \end{array} \end{equation}

	The future topology $\tp^+X $ is represented by its canonical local bases. In the last four 
columns we have applied the reflective and coreflective tools
    \begin{equation} \begin{array}{clr}
^\me\c \rho\Mtr \to \; !\rho\Mtr,   &   \; \pr_\infty\c \rho_\infty\Mtr \to \pSet  &   \text{(reflective),}
\\[5pt]
^\jo\c \rho\Mtr \to \; !\rho\Mtr,   &  \;\,  \pr_0\c \rho_\infty\Mtr \to \pSet   \q   &   \text{(coreflective).}
    \label{4.4.2} \end{array} \end{equation}

	We recall that the coreflective metric (and topology) of a $\rho$-space is always finer 
than the reflective one, as remarked in \ref{4.2}(b). Moreover $ x \prec_0 x' $ implies 
$ x' \prec_\infty x$.

\subsection{Reflective and coreflective analyses of maps}\label{4.5}
The Lipschitz maps between $\rho$-spaces can be studied using the previous tools.

\Ndt (a) {\em By reflective tools}. The $\rho_\infty$-maps $ f\c \de\bbR \to \de\bbR $ can be 
characterised as the mappings $ \bbR \to \bbR $ which are 

\LL

\ndt (i) Lipschitz (and uniformly continuous) for the euclidean metric $ d(x, x') = |x - x'|$,

\Ndt (ii) increasing for the natural order.

\LB

\skp	The necessity of (i) and (ii) is simply proved applying to $ f $ the `reflective tools' of 
\eqref{4.4.2} (while the coreflective ones would give no information). Conversely, an increasing 
map $ f\c \bbR \to \bbR $ which is $\la$-Lipschitz for the euclidean metric is also $\la$-Lipschitz 
for $\de\bbR$: for $ x \le x'$
    \begin{equation}
\la\de(x, x')  =  \la(x' - x)  =  \la|x - x'|  \; \ge \;  |fx - fx'|  =  \de(fx, fx'),
    \label{4.5.1} \end{equation}
while the case $ x' < x $ is trivially satisfied, since $ \la\de(x, x') = \infty $ (also for $ \la = 0$).

\Ndt (b) {\em By coreflective tools}. The $\rho_\infty$-maps $ f\c \de_0\bbR \to \de_0\bbR $ 
are characterised in the same way, using now the coreflective metric and coreflective preorder 
of $ \de_0\bbR$, which are the usual ones on the line.

	Again, we only have to check the sufficiency of conditions (i), (ii): an increasing map 
$ f\c \bbR \to \bbR $ which is $\la$-Lipschitz for the euclidean metric is also $\la$-Lipschitz for 
$ \de_0$: for $ x \le x'$
    \begin{equation}
\la\de_0(x, x')  =  \la(x' - x)  =  \la|x - x'|  \ge  |fx - fx'|  =  \de_0(fx, fx'),
    \label{4.5.2} \end{equation}
while the case $ x' < x $ gives $\la\de_0(x, x') = 0 = \de_0(fx, fx')$.

\Ndt(c) We have already seen that all Lipschitz maps $ \rho\bbR \to \rho\bbR $ are affine 
(in \ref{2.6}(a)), and the Lipschitz maps $ \rho\bbR \to \de_0\bbR $ are constant (in \ref{2.6}(b)).

\subsection{Another example}\label{4.6}
(a) A mapping $ f\c \de\bbR \to \rho\bbR $ is $\la$-Lipschitz if and only if
    \begin{equation}
f(x')  \le  f(x) + \la(x' - x),   \q   \for  x \le x',
    \label{4.6.1} \end{equation}
a form of {\em upper-right} uniform continuity condition: for every $ \ep > 0$, if $ x < x_0 + \ep $ 
then $ f(x) \le f(x_0) + \la\ep$. Thus, the function $f$ is $\la$-Lipschitz for the euclidean metric 
on each subinterval where it is increasing, but can have any `downward jump'.

\Ndt (b) For a mapping $ f\c \de\bbR \to \de_0\bbR $ the $\la$-Lipschitz condition still amounts 
to condition \eqref{4.6.1}.

\subsection{Lemma {\rm (Chaotic symmetrisations)}}\label{4.7}
{\em
(a) The reflective symmetrisation of an affordable non-positive $\rho$-metric space 
$ X $ is a chaotic $\rho$-space $ \rho_C|X|$.

\Ndt (b) In particular, $ (\rho\bbR)^\me = \rho_C\bbR, $ the chaotic $\rho$-line.
}
\begin{proof}
In fact, assuming that $ \rho(z, z') < 0 $ for a pair of points, we have for all $ x, x' $ in $ X$
$$
\hat{\rho}(x, x')  \le  \rho(x, z) + (2k - 1)\zeta + \rho(z', x'),
$$
where $ \zeta = \rho(z, z') \me \rho(z', z) < 0 $ and $ k $ is any integer $> 0$. (We are using the 
step-path $ \Ux = (x, z, z', ..., z, z', x')$, where the pair $(z, z')$ is repeated $k$ times.)
\end{proof}
%

\subsection{Linear metrics and topology}\label{4.8}
As remarked since the Introduction, real metrics have a complex relationship to topologies. 
We have so far considered four topologies on a $\rho$-metric space $X$, namely 
$ X^\jo \ge \tp^+X,  \tp^-X \ge X^\me$, but there are cases of interest where the `space' already 
bears a natural topology, even finer than $ X^\jo$.

	Extending the case of $ \rho\bbR$, let us consider a topological space $ X $ equipped 
with the linear $\rho$-metric $ \rho = \De\Phi $ associated to a {\em continuous} potential 
$ \Phi\c X \to \bbR$, as in \eqref{3.2.1}.

	The topology of $ X $ is (weakly) finer than the coreflective topology $ \tp^\jo(X, \rho)$, 
defined by the coreflective metric $ \check{\rho} = \rho \jo \rho\op$: in fact, every open ball of 
the latter
    \begin{equation}
B_\Phi(x_0, \ep)  =  \{x \in X  \sep  |\Phi(x) - \Phi(x_0)| < \ep\}   \qq   (\ep > 0),
    \label{4.8.1} \end{equation}
is also open in $ X$, as the $\Phi$-preimage of an open interval of $ \bbR$. 

	These topologies coincide for $\rho\bbR $ and many related cases, but the original 
topology can even be strictly finer than the coreflective one: for instance, if $ \Phi $ is constant, 
$ \rho = \check{\rho} = \hat{\rho} $ is the chaotic $\de$-metric, constant at $ 0$, independently 
of the topology of $ X$.

\section{Paths in real valued metric spaces}\label{s5}
As claimed in the Introduction, studying paths in a $\rho$-space $ X $ is likely the best way to 
explore its topological information and `direction'. We can now make precise many topics 
loosely described in the elementary example of \ref{0.2}.

	We are mostly interested in {\em Lipschitz paths}, that is, Lipschitz maps $ a\c \de\bbI \to X$, 
characterised by the condition $ ||a|| < \infty$. But we also define the {\em valuation} 
$ v(a) \le ||a|| $ of any mapping $ a\c \bbI \to X$, an extension of its length when $ X $ is a 
$\de$-space, studied in \cite{G3, G4}. Every Lipschitz path is thus {\em affordable}: 
$ v(a) < \infty$. The relationship between $\rho$-spaces and `spaces with directed paths' 
\cite{G1, G4} is dealt with in \ref{5.6}.

	$X $ always denotes a $\rho$-space.

\subsection{Lipschitz paths}\label{5.1}
(a) Paths (and homotopies) in $\rho$-spaces will be based on the standard $\de$-interval 
$ \de\bbI$. The positivity of its metric will play an important role.

	A {\em Lipschitz path} in a $\rho$-space $ X $ will be a Lipschitz map $ a\c \de\bbI \to X$. 
This means that there exists a (finite) $ \la \ge 0 $ such that
    \begin{equation}
\rho_X(a(t), a(t'))  \le  \la(t' - t),   \q   \for t \le t',
    \label{5.1.1} \end{equation}
where the case $ t > t' $ is automatically satisfied. We also say that $ a $ is a $\la$-{\em path}.

	The path $ a $ has a weight $ ||a|| \in [0, \infty[$, its least Lipschitz constant: see \ref{2.5}(b).

\Ndt (b) We note that

\LL

\ndt(i) the constant loop $ e_x $ at any point $ x \in X $ is Lipschitz, since $ \rho(x, x) $ is $ 0 $ 
or $ \Oinfty$,

\Ndt (ii) Lipschitz paths are closed under concatenation: concatenating two $\la$-paths gives 
a path admitting the Lipschitz constant $2\la$,

\Ndt (iii) Lipschitz paths are closed under (partial) reparametrisation, by any $\rho_\infty$-map 
$ \ph\c \de\bbI \to \de\bbI$; in other words, if $ a\c \de\bbI \to X $ is Lipschitz, so is obviously 
$ a\ph$.

\LB

\Ndt (c) A Lipschitz path $ \de\bbI \to \de\bbR $ is any increasing mapping which is Lipschitz for 
the euclidean metrics, as proved in \ref{4.5}(a) (for maps $ \de\bbR \to \de\bbR$). The same 
applies to the reparametrisations $\ph\c \de\bbI \to \de\bbI $ considered above. We also recall 
that a continuous endomap of $ \bbI $ need not be Lipschitz: the square root is not.

\Ndt (d) For a mapping $ a\c \de\bbI \to \rho\bbR $ the $\la$-Lipschitz condition amounts to
    \begin{equation}
a(x')  \; \le \;  a(x) + \la(x' - x),   \q\;\;  \for x \le x',
    \label{5.1.2} \end{equation}
and has been described in \ref{4.6} (for maps $ \de\bbR \to \rho\bbR$).

\Ndt (e) For a mapping $ a\c \de\bbI \to \de_0\bbR $ the $\la$-Lipschitz condition still amounts to 
condition \eqref{5.1.2}.

\subsection{The valuation of a path}\label{5.2}
More generally, we consider now a set-theoretical path $ a\c \bbI \to X $ (a mapping of sets) in 
a $\rho$-space. (A restriction to Lipschitz paths $ \de\bbI \to X $ would have no relevance here.)

\Ndt (a) We define the {\em valuation} $ v(a)$, also written as $ v_X(a) $ or $ v_\rho(a)$, by the 
following functions:
    \begin{equation}
v_\Ut(a)  =  \Sum_j \; \rho(a(t_{j-1}), a(t_j)),   \q   v(a)  =  \sup_\Ut \, v_\Ut(a)  \in  \Bw,
    \label{5.2.1} \end{equation}
where $\Ut = (t_j)$ stands for a finite increasing sequence  $ 0 = t_0 \le t_1 \le ... \le t_n =1$. 
The path $ a $ is said to be {\em affordable} if $ v(a) < \infty$. We shall see that every Lipschitz 
path is affordable (Theorem \ref{5.5}(f)).

	We always have $ v(a) \ge \rho(a(0), a(1))$. If $ X $ is a $\de$-metric space, $ v(a) \ge 0 $ 
is also called the {\em length} of $ a $ and written $ L(a)$, as in \cite{G3}, Definition \ref{1.2}, 
and in \cite{G4}, 6.1.8.

\Ndt (b) We define the {\em upper valuation}, or {\em total ascent}, of $a$ by its length in the 
$\de$-metric space $ X^+ $
    \begin{equation} \begin{array}{c}
v^+(a)  =  L_{X^+}(a)  =  \sup_\Ut \, \Sum_j \; \rho^+(a(t_{j-1}), a(t_j))  \in  \Bw^+,
\\[5pt]
v^+(a) \,  \ge  \, v(a).
    \label{5.2.2} \end{array} \end{equation}

	If $ v(a) < \infty$, we also define the {\em lower valuation}, or {\em total descent}, of $ a $ 
as the extended difference
    \begin{equation}
v^-(a)  =  v^+(a)  -  v(a)  \in  \Bw^+,
    \label{5.2.3} \end{equation}
which is $ \infty $ when $ v^+(a) = \infty $ or $ v(a) = \Oinfty $ (possibly both). Therefore
    \begin{equation}
v(a)  =  v^+(a)  -  v^-(a)   \qq  (\text{if } v(a) \jo v^+(a) < \infty).
    \label{5.2.4} \end{equation}

	If $ X $ is a $\de$-space, we are only interested in its length $ L(a) = v(a) = v^+(a) 
\in \Bw^+$. (Here $ v^-(a) = 0$, provided that $ v(a) < \infty$.)

\Ndt (c) {\em Absolute value}. If $ \rho $ is linear, its absolute value $ |\rho|(x, y) = |\rho(x, y)| $ 
is easily seen to be a finite symmetric $\de$-metric (generally called a pseudometric in the 
literature). 

	The {\em total variation} of the path $ a $ will be its $|\rho|$-valuation
    \begin{equation} \begin{array}{l}
L_{|\rho|}(a)  =  \sup_\Ux \, \Sum_j \, |\rho|(x_{j-1}, x_j)  =  \sup_\Ux \, \Sum_j \, |(\Phi(x_j) - \Phi(x_{j-1})|
\\[5pt]
\q\;\;\;\;  =  v^+_\rho(a) + v^-_\rho(a).
    \label{5.2.5} \end{array} \end{equation}
%
%

\subsection{Lemma}\label{5.3}
{\em
If $\rho$ is linear, with potential $\Phi$, for any path $a\c \bbI \to X$ (with reversed path $a\shp$)
}
    \begin{equation} \begin{array}{c}
v(a)  =  \rho(a(0), a(1))  =  \Phi(a(1)) - \Phi(a(0))  \in  \bbR,
\\[3pt]
v(a\shp)  =  - v(a),    \qq    v^-(a)  =  v^+(a\shp).
    \label{5.3.1} \end{array} \end{equation}
\begin{proof}
The first relation is obvious and the second is a consequence. As to to the third, we recall that 
$ \rho(y, x) = - \rho(x, y) $ and note that, for $\la \in \bbR$, we have $ \la = \la^+ - (- \la)^+$. 

	Letting $\Ut$ be any finite increasing sequence of points in $ \bbI$, with 
$0 = t_0 \le t_1 \le ... \le t_n =1$ and $p_j = a(t_j)$, we have by linearity the following relations:
    \begin{equation*} \begin{array}{c}
v(a)  =  \rho(p_0, p_1)  =  \Sum_j \, \rho(p_{j-1}, p_j)  
=  \Sum_j \, \rho^+(p_{j-1}, p_j) \! - \! \Sum_j \, \rho^+(p_j, p_{j-1}),
\\[5pt]
\Sum_j \, \rho^+(p_{j-1}, p_j)  =  v(a) + \Sum_j \, \rho^+(p_j, p_{j-1}),
\\[5pt]
v^+(a)  =  v(a) + \sup_\Ut \, (\Sum_j \, \rho^+(p_j, p_{j-1})))  =  v(a) + v^+(a\shp),
    \label{} \end{array} \end{equation*}
so that $ v^+(a\shp) $ is indeed $ v^+(a) - v(a) = v^-(a)$.
\end{proof}
%

\subsection{Examples and remarks}\label{5.4}
(a) For a path $ a $ in $ \de\bbR$, $L(a) = a(1) - a(0) $ if $ a $ is increasing, and 
$ L(a) = \infty $ otherwise.

\Ndt (b) For a path $ a $ in the linear $\rho$-space $ \rho\bbR$

\LL

\ndt - $v(a) = a(1) - a(0)$ is always finite, but $v^+(a) $ and $ v^-(a)$ can also be infinite 
(for instance, one can take the projection on a line of a Peano curve in the square),

\Ndt - also the weight $ ||a|| $ (of $ a\c \de\bbI \to \rho\bbR$) can be infinite: a classical example 
is $ a(t) = \sqrt{t}$,

\Ndt - if $a$ is increasing, $v^+(a) = v(a) = a(1) - a(0) $ and $ v^-(a) = 0$,

\Ndt - if $a$ is decreasing, $v^-(a) = - v(a) = a(0) - a(1) $ and $ v^+(a) = 0$.

\LB

\Ndt (c) For a path $ a $ in $ \de_0\bbR$, the length $L(a) $ is the upper valuation of $ a $ in 
$ \rho\bbR $ and can be infinite. For an increasing (resp.\ decreasing) path, $L(a) = a(1) - a(0)$ 
(resp.\ $ L(a) = 0$).

\Ndt (d) For a (set-theoretical) path $a$ in a discrete $\rho$-space $\rho_DX = \de_DX$, $L(a)$ 
is $ 0 $ or $ \infty$, according to the fact that $ a $ is constant or not.

	In a chaotic $\rho$-space $ \rho_CX$, $v(a) = \Oinfty$. In a chaotic $\de$-space 
$ \de_CX$, $L(a) = 0$.

\subsection{Theorem {\rm (Path valuation)}}\label{5.5}
{\em
Let $X$ be a $\rho$-metric space.

	The valuation of paths satisfies the following properties, where $a\c \bbI \to X$ is a 
set-theoretical path, $ e_x $ is the constant path at $ x \in X$, $a * b $ is the concatenation of 
two consecutive paths, $ \ph\c \bbI \to \bbI $ is a partial reparametrisation (an increasing map) 
and $ ||a|| \le \infty $ is the Lipschitz weight of $a\c \bbI \to X$ as a mapping defined on the 
$\de$-space $ \de\bbI$.
}
    \begin{equation*} \begin{array}{llr}
(a)   &   v(e_x)  =  \rho(x, x)  &   \text{(which is  0  {\em or }}  \Oinfty).
\\[4pt]
(b)   &   v(a * b)  =  v(a) + v(b).
\\[4pt]
(c)   &   v(a\ph)  \le  v(a),   &   \text{for a $\de$-space } X.
\\[4pt]
(d)   &   v(a) < \infty  \!  \implies \!  v(a\ph) < \infty,   &   \text{for a general $\rho$-space } X.
\\[4pt]
(e)   &   v(a\ph)  =  v(a),   &   \text{for an increasing {\em surjective} map } \ph.
\\[4pt]
(f)   &   v(a)  \le  ||a||   &   \text{(every Lipschitz path is affordable).}
\\[4pt]
(g)   &   v_Y(fa)  \le  \la v_X(a),   &   \text{for all $\la$-Lipschitz maps } f\c X \to Y.
    \label{} \end{array} \end{equation*}
\begin{proof}
There is a similar result in \cite{G4}, 6.1.8, for the positive case. Note that we can always replace 
a partition $ \Ut $ by a finer one $ \Ut'$, as the triangular property gives $ v_\Ut(a) \le v_{\Ut'}(a)$.

	Point (a) is obvious. For (b), the inequality $ v(a * b) \le v(a) + v(b) $ follows easily from the 
previous remark: given a partition $ \Ut $ for $ c = a * b$, call $ \Ut' $ its refinement by 
introducing the point $ t = 1/2 $ (if missing); then $ v_\Ut(c) \le v_{\Ut'}(c) \le v(a) + v(b)$. For the 
other inequality, it is sufficient to note that a partition for $ a $ and one for $ b $ yield a partition 
of $ [0, 2]$, which can be scaled down to the standard interval.

 	Point (e) follows from the fact that we are computing the mapping $ v $ on the same 
sequences: if $ (t_j) $ is an increasing sequence in $ \bbI$, $(\ph(t_j)) $ is also; 
conversely, given $ (t_j)$, there is an increasing sequence $ (s_j) $ such that 
$ \ph(s_j) = t_j$, for all indices $ j$. (Note also that $ \ph(0) = 0 $ and $ \ph(1) = 1$.)

	Points (c) and (d) are a consequence of (b) and (e). We can write $ a\si = a' * a\ph * a" $ 
for a suitable surjective reparametrisation $ \si$, so that $ v(a) = v(a') + v(a\ph) + v(a")$. In the 
positive case, all these lengths are $ \ge 0 $ and $ v(a) \ge v(a\ph)$. In the general case this 
can fail, but $ v(a\ph) = \infty $ still implies $ v(a) = \infty $ (since $ \infty $ is absorbent in 
$ \Bw$, even on $ \Oinfty$).

	Finally, points (f) and (g) are obvious, because
$$
\rho_X(a(t), a(t'))  \le  ||a||(t' - t),   \;\;   \rho_Y(fa(t), fa(t'))  \le  \la \rho_X(a(t), a(t')).
$$
\end{proof}
%

\subsection{Spaces with distinguished paths}\label{5.6}
This point is addressed to a reader with a basic knowledge of `spaces with directed paths', or 
{\em d-spaces}, our main structure for Directed Algebraic Topology, introduced in \cite{G1}, 
extensively studied in \cite{G4} and frequently used in the theory of concurrency 
(see \cite{FGHMR} and its references).

	In \cite{G3, G4} a $\de$-space $ X $ has an associated d-space ${\bf d} X $ on the 
same set, equipped with the reflective topology defined in \ref{4.1}, and directed paths consisting 
of all continuous mappings $ a\c \bbI \to X $ which are affordable, in the sense that $ v(a) < \infty$; 
in other words, their length $L(a)$ is finite.

	This can be extended to $\rho$-spaces (still using the reflective topology): the axioms of 
$\de$-spaces are still satisfied (by Theorem \ref{5.5}(a), (b), (d)), and $\rho$-maps preserve 
the directed paths (by \ref{5.5}(g)).

	Yet, we have already seen that this extension cannot be satisfactory: $\rho\bbR $ and several 
related $\rho$-spaces would get the chaotic topology. Using the final topology on $ X $ for the 
affordable paths might be preferable, but is not an extension of the positive case, and would 
give strange topologies, of little interest, to basic spaces like all $ \de\bbR^n$, for $ n > 1$.

\begin{small}

\skp	Essentially, this problem arises because we are working with `directed metrics': 
if $ X $ is a first-countable, locally path connected space, the final topology for {\em all} its 
(continuous) maps $ \bbI \to X $ is precisely the original topology, as proved in \cite{ChSW}, 
Proposition 3.11, in an article about diffeological spaces.

\end{small}

	Finally, this is another evidence of the vague relationship of real metrics to topology. 
As claimed in the Introduction, Part III will give a better solution to all these problems.

\section{Tensor products, other examples}\label{s6}
	The category $ \rho\Mtr $ has a `natural' symmetric monoidal closed structure, that extends 
the structure of $ \de\Mtr $ \cite{L1, G3, G4}, with the same formulas. We also examine the 
relationship between the tensor product and the cartesian product. The elementary example of 
the Introduction is reviewed in \ref{6.5}; the gravitational and elastic potential in \ref{6.6}, \ref{6.7}.

\subsection{Tensor product}\label{6.1}
(a) The tensor product $ X \te Y $ of $\rho$-spaces is the cartesian product of the underlying 
sets, endowed with the $l_1$-type $\rho$-metric (instead of the $l_\infty$-type $\rho$-metric of 
the categorical product)
    \begin{equation}
\rho((x, y), (x', y'))  =  \rho_X(x, x') + \rho_Y(y, y').
    \label{6.1.1} \end{equation}

	The unit is the $\de$-singleton $ \de\top $ (see \ref{2.4}(b)). On the other hand, 
$ X \te \rho\top = \rho_C|X|$, for all affordable $\rho$-spaces $ X$.

	This construction solves the usual universal problem, with respect to mappings which are 
1-Lipschitz in each variable. The exponential $ Z^Y $ is the set of 1-Lipschitz maps $ Y \to Z $ 
equipped with the $\rho$-metric
    \begin{equation}
\rho(h, k)  =  \sup_{y\in Y} \, \rho_Z(h(y), k(y))   \q\;\;   (h, k\c Y \to Z  \inn  \rho\Mtr).
    \label{6.1.2} \end{equation}

	The proof of the adjunction is standard, and can be found in point (d). The tensor powers 
of the basic $\rho$-structures of the line are written as $ \de\bbR^{\te n}$, $\rho\bbR^{\te n}$ 
and $ \de_0\bbR^{\te n}$.

\Ndt (b) {\em Restricting to positive } $\de$-{\em metrics}, the cartesian and tensor product 
satisfy the following inequalities:
    \begin{equation}
X \ti Y  \le  X \te Y  \le  2(X \ti Y)   \qq   (\text{in } \de\Mtr).
    \label{6.1.3} \end{equation}

	In $ \de_\infty\Mtr $ these products are thus isomorphic, and denote isomorphic functors in 
two variables. This is not true in $ \de\Mtr$, and we shall always distinguish such objects (and 
functors): the notation $ X \ti Y $ (resp.\ $ X \te Y$) will always denote the structure given by the 
$l_\infty$-type (resp.\ $l_1$-type) $\de$-metric.

\Ndt (c) For $\rho$-metrics we still have
    \begin{equation}
X \te Y  \le  2(X \ti Y)   \qq   (\text{in }\rho\Mtr),
    \label{6.1.4} \end{equation}
but the first inequality in \eqref{6.1.3} fails (also for symmetric $\rho$-spaces). For instance, the 
$\rho$-spaces $ \de\top \ti \rho\top = \de\top $ and $ \de\top \te \rho\top = \rho\top $ give the 
reverse inequality, and are not even Lipschitz-isomorphic.

\Ndt (d) Finally, to prove the adjunction, we verify the bijective correspondence 
$ f \leftrightarrow g $ between 1-Lipschitz maps
    \begin{equation}
f\c X \te Y \to Z,   \q   g\c X \to Z^Y   \qq   (\text{in } \rho\Mtr),
    \label{6.1.5} \end{equation}
linked by the relation $ f(x, y) = g(x)(y)$, for all $ x \in X$, $y \in Y$. (These quantifiers are left 
understood, below).

	Given $f$, the relation $\rho_Z(f(x, y), f(x', y')) \le \rho_X(x, x') + \rho_Y(y, y')$ 
shows that

\LL

\ndt -  every map $ g(x) = f(x, -)\c Y \to Z $ is 1-Lipschitz, because (taking into account that 
$ \rho_X(x, x) \le 0$)
$$
\rho_Z(f(x, y), f(x, y'))  \le  \rho_X(x, x) + \rho_Y(y, y')  \le  \rho_Y(y, y'),
$$

\ndt - the map $ g\c X \to Z^Y $ is 1-Lipschitz, because we have, in the same way
    \begin{equation*} \begin{array}{c}
\rho_Z(g(x)(y), g(x')(y))  =  \rho_Z(f(x, y), f(x', y))
\\[3pt]
\qqq    \le \;  \rho_X(x, x') + \rho_Y(y, y)  \le  \rho_X(x, x').
    \label{} \end{array} \end{equation*}

\LB

	The other implication is similarly proved, inserting a point $ (x', y) $ between $ (x, y) $ and 
$ (x', y')$.

\subsection{Theorem {\rm (Properties of the tensor product)}}\label{6.2}
{\em
Let $ X, Y $ be $\rho$-metric spaces.

\Ndt (a) The reflector $ ^\me\c \rho\Mtr \to \; !\rho\Mtr $ is a strict monoidal functor
    \begin{equation}
(X \te Y)^\me  =  X^\me \te Y^\me,
    \label{6.2.1} \end{equation}
and forms, with the embedding $ !\rho\Mtr \sub \rho\Mtr$, a strict adjunction of monoidal 
categories.

\Ndt (b) The coreflector $ ^\jo\c \rho\Mtr \to \; !\rho\Mtr $ is a lax monoidal functor
    \begin{equation}
(X \te Y)^\jo  \le  X^\jo \te Y^\jo,
    \label{6.2.2} \end{equation}
with comparison $ X^\jo \te Y^\jo \to (X \te Y)^\jo$. It forms a strict/lax adjunction.

\Ndt (c) As to the reflective and coreflective preorders (in \ref{4.3})
    \begin{equation}
(x \prec_\infty x'  \inn  X, \;  y\prec_\infty y'  \inn  Y)  \iff   (x, y) \prec_\infty (x', y')   \inn  X \te Y,
    \label{6.2.3} \end{equation}
    \begin{equation}
(x \prec_0 x'  \inn  X, \;   y\prec_0 y'  \inn  Y)	\implies	(x, y) \prec_0 (x', y')   \inn  X \te Y.
    \label{6.2.4} \end{equation}

\Ndt (d) For a set-theoretical path $ \lan a, b\ran \c \bbI \to X \te Y$
}
    \begin{equation}
v(\lan a, b\ran)  =  v(a) + v(b).
    \label{6.2.5} \end{equation}
\begin{proof}
We replace the pair $ X, Y $ with a finite family of $\rho$-metric spaces $ X_1,..., X_n$; this will 
give a simpler notation and a clearer proof. An element of the set $ \Pro_i |X_i| $ is written as 
$ \Ux = (x_i)$. We use the notation $ \rho_X(x, x') = X(x, x')$, which comes from viewing a 
$\rho$-metric space as an enriched category (in \ref{2.1}).

\Ndt (a) We now prove that $ (\Ten_i \, X_i)^\me = \; \Ten_i \, X_i^\me$. The $\rho$-metric of the 
latter satisfies
$$
\Ten \, X_i^\me (\Ux, \Uy) = \Sum_i \, X_i^\me (x_i, y_i)  \; \le \; \Sum_i \, X_i(x_i, y_i) 
=  (\Ten \, X_i)(\Ux, \Uy).
$$

	Since $\, \Ten \, X_i^\me $ is symmetric, we have $\, \Ten \, X_i^\me \le (\Ten \, X_i)^\me$. 
The reverse inequality is more subtle: take a sequence of $ n+1 $ points $ \Uz^j$, which 
varies from $ \Ux $ to $ \Uy $ by changing one coordinate at a time
$$
\Uz^j  =  (y_1,..., y_j, x_{j+1},..., x_n),   \q   \Uz^0  =  \Ux, \;\,  \Uz^n  =  \Uy   \q   (j = 0,..., n).
$$
and apply the triangular inequality
$$
(\Ten \, X_i)^\me(\Ux, \Uy) \;  \le \; \Sum_j \, (\Ten \, X_i)^\me(\Uz^{j-1}, \Uz^j).
$$

	Now, $ \Uz^{j-1}$ and $ \Uz^j $ only differ at the $j$-th coordinate ($x_j $ or $ y_j$, 
respectively); restricting the domain of the `inf' in the right term above to those sequences in 
$ (\Ten \, X_i)^\me $ where only the $j$-th coordinate changes, we get the $\de$-metric 
$ X_j^\me$, and the inequality
$$
\Sum_j \, (\Ten \, X_i)^\me(\Uz^{j-1}, \Uz^j)  \; \le \; \Sum_j \, X_j^\me(x_j, y_j)  =  
(\Ten_j \, X_j)^\me(\Ux, \Uy).
$$

\Ndt(b) In fact
    \begin{equation*} \begin{array}{c}
(\Ten \, X_i)^\jo(\Ux, \Uy)  =  \Sum_i \, X_i(x_i, y_i)) \jo (\Sum_i \, X_i(y_i, x_i)
\\[5pt]
\qq\q   \le \;  \Sum_i \, (X_i(x_i, y_i) \jo X_i(y_i, x_i))  =  \Ten \, X_i^\jo(\Ux, \Uy).
    \label{} \end{array} \end{equation*}

	Finally, point (c) is obvious. For \eqref{6.2.5}, the inequality $ v\lan a, b\ran \le v(a) + v(b) $ 
is obvious, and the other follows from the first remark in the proof of Theorem \ref{5.5}: given a 
partition for $ a $ and one for $ b$, by using a common refinement for both we get higher values.
\end{proof}

\subsection{Corollary I {\rm (Reflective symmetrisation and products)}}\label{6.3}
{\em
(a) For a finite family of $\rho$-metric spaces $ X_1,..., X_n$
    \begin{equation}
\Pro X_i^\me  \le  (\Pro X_i)^\me,  \q  (\Ten X_i)^\me  =  \Ten X_i^\me  \le \; n(\Pro X_i^\me).
    \label{6.3.1} \end{equation}

\ndt (b) More particularly, for $\de$-spaces
    \begin{equation}
\Pro X_i^\me  \le  (\Pro X_i)^\me  \le  \, (\Ten X_i)^\me  =  \Ten X_i^\me  \le  n(\Pro X_i^\me).
    \label{6.3.2} \end{equation}

	All these symmetric $\de$-metrics are Lipschitz-equivalent and induce the same topology: 
the functors $ ^\me\c \de_\infty\Mtr \to \; !\de_\infty\Mtr $ and $ ^\me\c \de_\infty\Mtr \to \Top $ 
preserve finite products and take tensor products to cartesian products.
}
\begin{proof}
(a) To compare $ \Pro X_i^\me$ and $ (\Pro X_i)^\me$, we note that
$$
(\Pro X_i^\me)(\Ux, \Uy)  =  \sup_i \, X_i^\me (x_i, y_i)  \le  \sup_i \, X_i(x_i, y_i)  
  =  (\Pro X_i)(\Ux, \Uy).
$$

	Since $ \Pro(X_i^\me) $ is symmetric, it follows that $ \Pro(X_i^\me) \le (\Pro X_i)^\me$. 
The second part follows from \eqref{6.2.1}.

\Ndt (b) For $\de$-spaces, the inequality $ (\Pro X_i)^\jo \le (\Ten X_i)^\jo $ follows from 
$ \Pro X_i \le \Ten X_i$, proved in \eqref{6.1.3}.
\end{proof}
%

\subsection{Corollary II {\rm (Coreflective symmetrisation and products)}}\label{6.4}
{\em
(a) For a finite family of $\rho$-spaces $X_1,..., X_n$
    \begin{equation}
\Pro X_i^\jo  =  (\Pro X_i)^\jo,   \q
	(\Ten X_i)^\jo  \le  \Ten X_i^\jo  \le  n(\Pro X_i^\jo).
    \label{6.4.1} \end{equation}

	Therefore the functor $ ^\jo\c \rho_\infty\Mtr \to \Top $ preserves finite products (but takes 
arbitrary products to the box topology).

\Ndt (b) More particularly, for $\de$-spaces
    \begin{equation}
\Pro X_i^\jo  =  (\Pro X_i)^\jo  \le  (\Ten X_i)^\jo  \le  \Ten (X_i^\jo)  \le  n(\Pro X_i^\jo).
    \label{6.4.2} \end{equation}

	Again, these symmetric $\de$-metrics are Lipschitz-equivalent and induce the same 
topology: the functor $ ^\jo\c \de_\infty\Mtr \to \Top $ preserves finite products, and takes tensor 
products to cartesian products.

\Ndt (c) If the symmetrisations $ X_i^\jo $ of the $\rho$-spaces $ X_i $ are {\em positive}
    \begin{equation}
\Pro X_i^\jo  =  (\Pro X_i)^\jo  \le  \Ten (X_i^\jo)  \le  n(\Pro X_i^\jo),
    \label{6.4.3} \end{equation}
and again these $\rho$-spaces are Lipschitz-equivalent (but $(\Ten X_i)^\jo $ need not be).
}
\begin{proof}
(a) Coreflective symmetrisation is a right adjoint, and preserves products. The relation 
$ (\Ten X_i)^\jo \le \Ten X_i^\jo $ is proved in \ref{6.2}, and the last inequality in \eqref{6.4.1} 
is obvious, as in \eqref{6.1.4}.

\Ndt (b) As in Corollary I, the inequality $ (\Pro X_i)^\jo \le (\Ten X_i)^\jo $ follows from
 $ \Pro X_i \le \Ten X_i$, proved in \eqref{6.1.3} for $\de$-spaces. We know that this does 
 not extend to $\rho$-spaces, even in the symmetric case.

\Ndt (c) If all $ X_i^\jo $ are positive, the inequality $ \Pro X_i^\jo \le \Ten (X_i^\jo) $ is a 
straightforward consequence of \eqref{6.1.3}, and we already know that 
$ (\Pro X_i)^\jo = (\Pro(X_i^\jo)$.
\end{proof}
%

\subsection{Ascent and descent}\label{6.5}
(a) Re-examining the elementary example of the Introduction, in \ref{0.2}, we consider the 
$\rho$-space $ X = \de\bbR \ti \rho\bbR$, with the product $\rho$-metric
$\rho((x, y), (x', y'))  =  \de(x, x') \jo \rho(y, y')$
   \begin{equation} 
\rho((x, y), (x', y')) \;  =  \;
    \begin{cases}
(x' - x) \jo (y' - y),  &  \IF         x \le x',
\\[3pt]
\infty                      &  otherwise.
    \end{cases}
    \label{6.5.1} \end{equation}   

\Ndt (b) A $\rho_\infty$-map $ a = \lan a_1, a_2\ran \c \de\bbI \to X$
%
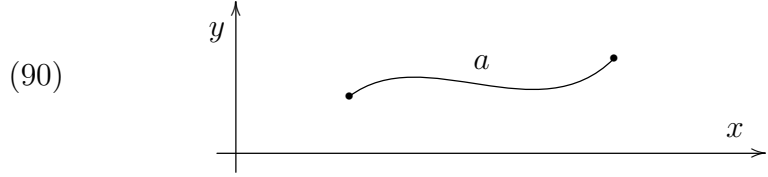
    \begin{equation} 
\xy <.5mm, 0mm>:
(45,4) *{a}; (112,-14) *{x}; (-25,12) *{y}; 
(0,28) *{}; (0,-32) *{}; 
@i @={(10, -5), (80, 5)} @@{*{\bu}};
(10, -5); (80, 5) **\crv{(30,10)&(60,-15)},
\POS(-25,-20) \ar+(145,0),  \POS(-20,-25) \ar+(0,45),  
\endxy
    \label{6.5.2} \end{equation}
has components $ a_1, a_2 $ such that

\LL

\ndt - the path $ a_1\c \de\bbI \to \de\bbR $ is any increasing map $ \bbI \to \bbR $ which 
is Lipschitz for the euclidean metric, as proved in \ref{4.5}(a),

\Ndt - the path $ a_2\c \de\bbI \to \rho\bbR $ satisfies the Lipschitz condition, examined in 
\ref{4.6} and \ref{5.1}(d).

\LB

\Ndt (c) Evaluation, total ascent, total descent and total variation have already been examined 
in \ref{0.2}, consistently with the definitions of \ref{5.2}.

\subsection{Gravitational potential and metric}\label{6.6}
(a) The gravitational potential around the Earth, for a unit mass, can be expressed as
    \begin{equation}
V(r)  =  - k r_0/r,    \qq   \for  r \ge r_0,
    \label{6.6.1} \end{equation}
where $ r_0 $ is radius of the Earth, $ r \ge r_0 $ is the distance from the centre of the Earth, and 
$ k = GM/r_0 $ is a positive constant; $ G $ is the gravitational constant and $ M $ is the mass of 
the Earth.

	We rewrite it in the form
    \begin{equation}
\Phi(x)  =  k - k/(1 + x)  =  kx/(1 + x),    \q\;\;   \for  x \ge 0,
    \label{6.6.2} \end{equation}
where $ x = r/r_0 - 1 \ge 0 $ is the altitude, measured with respect to the Earth radius; we have 
added the constant $ k $ so that the potential $ \Phi $ annihilates at the Earth surface: 
$ \Phi(0) = 0$. The function $ \Phi $ is increasing.

	The associated (linear) $\rho$-metric $\rho(x, y)  =  \Phi(y) - \Phi(x)$
measures the gain of potential energy from altitude $ x $ to altitude $ y$, or the opposite 
of the work of the gravitational field in this transition.

\Ndt (b) At small altitude, that is, for small $x$, one can use a linear approximation $\psi(x) = kx$ 
(which amounts to taking the scalar gravitational field constant), getting a $\rho$-metric
    \begin{equation}
\rho'(x, y)  =  \psi(y) - \psi(x)  =  k(y - x)    \q\;\;   \for  x, y \ge 0,
    \label{6.6.3} \end{equation}
which is linear and invariant up to (positive) translations.

\Ndt (c) Consider now a spherical mass $ M $ of radius $ r_0$, centred at the origin of $ \bbR^3$. 
Its gravitational potential, for a unit mass, can be expressed as
    \begin{equation}
V(\Ux)  =  - kr_0/||\Ux||,    \qq   \for  ||\Ux|| \ge r_0.
    \label{6.6.4} \end{equation}

	An interested reader can extend this to the gravitational (or electric) potential generated by 
a finite distribution of spherical masses (or charges), in the euclidean space.

\subsection{Elastic potential and metric}\label{6.7}
(a) In the cartesian plane $ xy $ a (linear) elastic field centred at the origin gives a potential function
    \begin{equation}
\Phi(x, y)  =  \la(x^2 + y^2),   \q   \for  (x, y) \in \bbR^2,
    \label{6.7.1} \end{equation}
where $ k = 2\la > 0 $ is the elastic constant of the field: $ F(x, y) = - k(x, y)$.
The associated (linear) $\rho$-metric
    \begin{equation}
\rho((x, y), (x', y'))  =  \Phi(x', y') - \Phi(x, y)  =  \la((x'^2 + y'^2) - (x^2 + y^2)),
    \label{6.7.2} \end{equation}
measures the gain of potential energy from position $(x, y)$ to $(x', y')$, that is the 
opposite of the work of the elastic field in this transition.

	Letting $\rho_E\bbR $ have $ \rho_E(x, x') = x'^2 - x^2$, we are considering the linear 
$\rho$-space $ (\la\rho_E\bbR)^{\te 2}$.



\begin{thebibliography}{99}

\bibitem[Ba]{Ba} M. Barr, *-Autonomous categories, Lecture Notes in Mathematics, Vol. 752, 
Springer, 1979.

\bibitem[Bk]{Bk} N. Bourbaki, Topologie g\'en\'erale, Ch. 10, Hermann, 1961.

\bibitem[ChSW]{ChSW} J.D. Christensen, G. Sinnamon and E. Wu, The D-topology for diffeological 
spaces, Pacific. J. Math. 272 (2014), 87--110.

\bibitem[EK]{EK} S. Eilenberg and G.M. Kelly, Closed categories, in Proc. Conf. Categorical 
Algebra, La Jolla 1965, Springer, 1966, pp. 421--562. 

\bibitem[FGHMR]{FGHMR} L. Fajstrup, E. Goubault, E. Haucourt, S. Mimram and M. Raussen, 
Directed algebraic topology and concurrency, Springer, 2016.

\bibitem[Fl]{Fl} R. Flagg, Quantales and continuity spaces, Algebra Universalis, 37 (1997), 
257--276.

\bibitem[G1]{G1} M. Grandis, Directed homotopy theory, I. The fundamental category, 
Cah. Topol. G\'eom. Diff\'er. Cat\'eg. 44 (2003), 281--316. Available at

https://www.numdam.org/issues/CTGDC\_2003\_44\_4/

\bibitem[G2]{G2} M. Grandis, Categories, norms and weights, J. Homotopy Relat. Struct. 
2 (2007), No. 9, 171--186.

\bibitem[G3]{G3} M. Grandis, The fundamental weighted category of a weighted space 
(From directed to weighted algebraic topology), Homology Homotopy Appl. 9 (2007), 
221--256. Available at

https://intlpress.com/JDetail/1805807194991902721

\bibitem[G4]{G4}  M. Grandis, Directed Algebraic Topology, Models of non-reversible worlds, 
Cambridge Univ. Press, 2009. Available at

https://www.researchgate.net/publication/267089582

\bibitem[G5]{G5}  M. Grandis, Weighted algebraic topology, II (Relative metrics), 
Preprint, arXiv, 2026.  

\bibitem[Ke]{Ke} J.C. Kelly, Bitopological spaces, Proc. London Math. Soc. 13 (1963), 
71--89.

\bibitem[Ky]{Ky} G.M. Kelly, Basic concepts of enriched category theory, Cambridge 
University Press, 1982.

\bibitem[L1]{L1} F.W. Lawvere, Metric spaces, generalized logic and closed categories, 
Rend. Sem. Mat. Fis. Univ. Milano 43 (1974), 135--166. Republished in: Reprints 
Th. Appl. Categ. 1 (2002), 1--37. Available at

http://www.tac.mta.ca/tac/reprints/articles/1/tr1.pdf

\bibitem[L2]{L2} F.W. Lawvere, State Categories, Closed Categories, and the Existence 
Semi-Continuous Entropy Functions, IMA Preprint Series 84, University of Minnesota, 1984.

\bibitem[Ma]{Ma} S. Mac Lane, Categories for the working mathematician, Springer, 1971.

\bibitem[Mu]{Mu} C.J. Mulvey, ``\&", Rend. Circ. Mat. Palermo, Suppl. 12 (1986), 99--104

\bibitem[Ro]{Ro} K. Rosenthal, Quantales and Their Applications, Longman Scientific \& 
Technical, 1990.

\bibitem[Wi]{Wi} S. Willerton, The Legendre-Fenchel transform from a category theoretic 
perspective, arXiv 1501.03791v1, 2015.


\bibitem[Ye]{Ye} D. Yetter, Quantales and (noncommutative) linear logic, J. Symb. Log. 55 
(1990), 41--64.

\end{thebibliography}
\end{document}